

\baselineskip=14pt
\parskip=10pt

\font\eightrm=cmr8 

\magnification=\magstephalf

\def\P{{\cal P}}

\def\1{{\overline{1}}}
\def\2{{\overline{2}}}
\def\d{{\overline{d}}}
\def\k{{\overline{k}}}
\parindent=0pt
\overfullrule=0in

\def\frac#1#2{{#1 \over #2}}
\centerline
{\bf  
A Multi-Computational Exploration of Some Games of Pure Chance
}
\bigskip
\centerline
{\it Thotsaporn ``Aek'' THANATIPANONDA and Doron ZEILBERGER}
\bigskip

{\bf Abstract.}
In the spirit of ``multi-culturalism'', we use four kinds of computations: simulation, numeric, symbolic, 
and ``conceptual'', to explore some ``games of pure chance'' inspired by children board games like ``Snakes and Ladders'' 
(aka ``Chutes and Ladders'') and ``gambler's ruin with unlimited credit''. Even more interesting than the many computer-generated 
actual results described in this paper and its web-site extension, is our broad-minded, ecunemical approach, not favoring, 
a priori, any one of the above four kinds of computation, but showing that, a posteriori, symbolic computation is the most 
important one, since (except for simulation) numerics can be made more efficient with the help of symbolics (in the ``downward'' direction), 
and, (in the ``upward'' direction) the mere existence of certain symbolic-computational algorithms imply interesting ``qualitative'' results, 
that certain numbers are always rational, or always algebraic, and certain sequences are always polynomial, or C-recursive, 
or algebraic, or holonomic. This article is accompanied by four Maple packages, and numerous input and output files, 
that readers can use as templates for their own investigations. 

{\bf The Maple packages.} This article is accompanied by four Maple packages 

$\bullet$ {\tt SnakesAndLadders.txt}  \quad ;

$\bullet$ {\tt PosPileGames.txt} \quad ;

$\bullet$ {\tt GenPileGames.txt} \quad ;

$\bullet$ {\tt VGPileGames.txt} \quad .

They are available, along with numerous input and output files, from the front of this article

{\tt http://www.math.rutgers.edu/\~{}zeilberg/mamarim/mamarimhtml/chance.html} \quad .

{\bf Chapter 0: Important Background and Definitions} 

{\bf The three Kinds of Games}

First we have TicTacToe, Checkers, Chess, and Go, that are games of {\bf no} chance, i.e. you don't need
Lady Luck to help you (except the luck to be born smart).

Then we have  games of {\bf some} chance, combining luck and skill, like Poker, Bridge, Backgammon, and countless other games.

Finally, we have games of {\bf pure} chance, where there is no skill involved,  like the children board games ``Snakes and Ladders'',
``Candy Land'', and gambling in a casino (where at each round you always bet the {\bf same} amount, independent of your capital).
Also the final stage of Backgammon, where all the pieces are in their final quarter, may be considered
as a game of pure chance (assuming that both players follow a fixed bear-off strategy, like the greedy algorithm).

In this article we consider certain families of games of pure chance, to be defined later.

{\bf The four Kinds of Computation to study Games of Pure Chance }

$\bullet$ 
First we have  {\bf simulation} (aka Monte Carlo), where you let the computer play the game many times,
and to estimate the probability of winning you divide the number of games won by the total number of games,
and estimate the expected length of the game by the sample average. This is very unreliable, especially,
if the variance is large. Nevertheless, it is a useful check of the more accurate results
obtained by other kinds of computations.

$\bullet$ Next we have {\bf numeric} computation, whose output is a {\bf number}. For
example the probability that I will get $10$ Heads if I toss a fair coin $20$ times is the
{\bf exact number} $\frac{46189}{262144}$ that equals  $0.176197052\dots$. This is a {\bf numerical} answer
computed by a numeric computer programming language by plugging in the parameters $n=20$ and $k=10$,
after a human coder {\bf hard-wired} the human-generated (by human geniuses Pascal, Fermat,
and possibly 13th century Chu) formula $\frac{n!}{2^n k!(n-k)!}$.

$\bullet$ Next we have {\bf symbolic} computation, where the computer generates general formulas, that
are valid for {\bf symbolic} parameters. Here the human coder designs general purpose
algorithms that can output many such formulas (or more general schemes) that incorporate
`infinitely many facts'. For example nowadays the formula   $\frac{n!}{2^n k!(n-k)!}$ can be gotten
{\it ab initio} using the Almkvist-Zeilberger algorithm [AZ]. Just type

{\tt AZd((1+x)**n/x**(k+1)/x,x,n,N)[1];} \quad ,

in the Maple package {\tt http://sites.math.rutgers.edu/\~{}zeilberg/tokhniot/EKHAD.txt }.

Note that in some sense, numeric computation is a special case of symbolic computation, since numbers {\it are} symbols.

$\bullet$ Finally we have {\bf conceptual computation}, that sounds like an {\it oxymoron}.
In this case, we don't compute {\bf exact numbers}, or {\bf explicit symbolic expressions}, but nevertheless
use (often mental) ``meta-computations'' to assert the {\it nature} of the desired number.
Is it rational? Is it algebraic? Or the {\it nature} of the desired sequence.
Is our sequence polynomial?, or C-recursive?, or   P-recursive?, or ``none of the above''?

Before going on, let's remind us what these mean.

{\bf The Four kinds of Numbers  }

$\bullet$ First we have {\bf positive integers}, aka  {\it natural numbers} (and zero), that were created by God.

The socially-constructed numbers are as follows.

$\bullet$ {\bf  Positive rational numbers}: that are of the form $\frac{m}{n}$ where $m$ and $n$ are positive integers and $n \neq 0$.
Pythagoras believed that all numbers are rational.

$\bullet$ {\bf Algebraic numbers}: numbers $x$ that satisfy a polynomial equation $P(x)=0$, where $P$ is
a polynomial with {\bf integer} coefficients. For example $\sqrt{3}$ and $1+ \sqrt{5}$, that we will
encounter later in this article, are algebraic numbers, since they are roots of the equations
$x^2-3=0$ and $x^2-2x-4=0$ respectively.

$\bullet$ {\bf Computable numbers}, numbers $x$, such that for any $\epsilon>0$, you can compute (hopefully fast)
a rational number $x_0$ such that $|x-x_0| \leq \epsilon$.

We have no use for any other numbers. Note that the above families are `strictly increasing'.
Every integer is rational (take $n=1$). Every rational number is algebraic (take the degree of $P$ to be $1$).
Every algebraic number is computable (using approximation algorithms like Newton-Raphson, or Horner).

We will also talk about {\bf sequences} of numbers.

{\bf Four Kinds of Sequences}

$\bullet$ Polynomial sequences, $a(n)=p(n)$, where $p(n)$ is a polynomial in $n$. For example $\{n\}$, $\{n^2\}$, $\{n^{1000}\}$.

$\bullet$ $C-$recursive (aka $C$-finite) sequences, $a(n)$, that satisfy  a {\bf linear recurrence equation with constant coefficients}
$$
a(n)=c_1 a(n-1)+ \dots + c_L a(n-L) \quad,
$$
for some positive integer $L$ and some constants $c_1, \dots, c_L$. Equivalently, the (ordinary) {\bf generating function}
$$
\sum_{n=0}^{\infty} a(n) t^n \quad,
$$
is a {\bf rational function} of $t$, i.e. can be written as $P(t)/Q(t)$ for some {\bf polynomials} $P(t),Q(t)$.
The most famous $C$-finite sequences (that are not polynomials) are $\{2^n\}$ and the sequence of Fibonacci numbers $\{F_n\}$.

See [Z1], [Z2] and [KP] about them.
Note that any polynomial sequence is $C$-recursive. The denominator of its generating function is $(1-t)^{d+1}$, where $d$ is the degree.

$\bullet$  Algebraic sequences that satisfy a {\bf non-linear} recurrence with constant coefficients. 
Equivalently sequences whose ordinary generating function
$$
f(t) := \, \sum_{n=0}^{\infty} a(n) t^n \quad,
$$
satisfy an equation of the form
$$
P(t,f(t))=0 \quad,
$$
for some {\bf polynomial} of {\bf two} variables $P(x,y)$.
The most famous algebraic sequence (that is not $C$-recursive) is the sequence of {\bf Catalan numbers} $\frac{(2n)!}{n!(n+1)!}$, 
whose generating function satisfies $f(t)=1+tf(t)^2$. See [KP] chapter 6, and [Z2].
Every $C$-recursive sequence is algebraic, where the defining equation, $P(t,f(t))=0$, has degree one in $f(t)$.

$\bullet$ $P-$recursive (aka {\it [discrete] holonomic}) sequences, $a(n)$, that satisfy  a {\bf linear recurrence equation with polynomial coefficients}
$$
c_0(n)\,a(n) + c_1(n) a(n-1) + \dots + c_L(n) a(n-L) \, = \, 0\quad,
$$
for some positive integer $L$ and some {\bf polynomials}, in the discrete variable $n$,
$c_0(n), \dots, c_L(n)$. Equivalently, the (ordinary) {\bf generating function}
$$
\sum_{n=0}^{\infty} a(n) t^n \quad,
$$
is $D$-finite (aka [continuous] holonomic), i.e. it satisfies a {\bf linear differential equation with polynomial coefficients}:
$$
d_0(t)\,f(t) + d_1(t) f'(t) + \dots + d_M(t) f^{(M)}(t) \, = \, 0\quad,
$$
for some positive integer $M$ and some {\bf polynomials} $d_0(t), \dots, d_M(t)$. 
Thanks to a famous theorem (see [KP], Theorem 6.1), every  algebraic sequence is also $P$-recursive.
Note that the converse is {\bf not} true, e.g. $\{\frac{(3n)!}{n!^3}\}$.

We need to define one more kind of {\bf number}.

{\bf Definition}: A real number is a {\bf holonomic constant} if it is the sum of a convergent  power series
$$
\sum_{n=0}^{\infty} \,  a_n \quad ,
$$
where the sequence $a_n$ satisfies  a {\bf linear recurrence equation with polynomial coefficients} (in $n$) whose coefficients are
polynomials (in $n$) with {\bf integer coefficients}, and the initial conditions are {\bf rational numbers}.
Note that, just like rational numbers and algebraic numbers, being a holnomic constant is a {\it big deal}, since
there are only a {\bf countable} number of them.

{\bf Chapter 1: Specific Games of Pure Chance inspired by the board game ``Snakes and Ladders''}

The game of ``Snakes and Ladders'' consists of $100$ squares, and players take turns spinning a spinner
(or equivalently rolling a die) with six equally likely outcomes $1,2,3,4,5,6$. If the player is currently
at location $i$, and got $j$, then she goes to location $i+j$. In addition there are several
``snakes'' (chutes) $[i,j]$ where $i>j$ and if the player landed on location $i$ he must jump down to location $j$.
There are also several ``ladders'' $[i,j]$ with $i<j$, where the player jumps forward from location $i$ to location $j$.

To be specific, the version called ``Chutes and Ladders''  manufactured by {\it Winning Moves Games}, has
$100$ locations, with the following $10$ chutes
$$
[16, 6] \,, \, [47, 26] \, , \, [49, 11] \, , \, [56, 53] \, , \, [62, 19] \, , \, [64, 60] \, , \, [87, 24] \, , \, [93, 73]
\, , \, [95, 75] \, , \, [98, 78] \quad,
$$
and the following $9$ ladders:
$$
[1, 38] \, , \, [4, 14] \, , \, [9, 31] \, , \, [21, 42] \, , \, [28, 84] \, , \, [36,44] \, , \, [51, 67] \, , \, [71, 91] \, , \, [80, 100] \quad .
$$
On the other hand the version called ``Snakes and Ladders'' manufactured by {\it Cardinal Industries} has the
following $7$ chutes (that are called there ``snakes'')
$$
[16, 5] \, , \, [50, 8] \, , \, [63, 20] \, , \, [57, 25] \, , \, [98, 46] \,, \, [98, 76] \, , \, [95, 90] \quad ,
$$
and the following $8$ ladders
$$
[2, 44] \, , \, [6, 13] \, , \, [9, 31] \, , \, [28, 84] \, , \, [59, 61] \, , \, [67, 93] \, , \, [70, 73] \, , \, [79, 100] \quad .
$$

All such games can be modeled in terms of a (finite) Markov process with one {\it absorbing state}.

We have a {\it directed graph} with  a set of vertices $V$, one of whose vertices is the
absorbing state, $e$. There are no edges out of $e$, and out of each vertex $v \in V$ there is a set
of outgoing neighbors, let's call it $N(v)$, so that there is a directed edge from $v$ to
each member $u \in V$, and a probability distribution on $N(v)$, $\{ p_{vu} \, | \, u \in N(v)\}$, such that
$\sum_{ u \in N(v)} p_{vu}=1$.

Since this is a game of {\bf pure} chance, where there is no strategy, and it is essentially
a race between the players, rather than having the huge
``state'' space $V \times V$, it is much more efficient to first consider the {\it solitaire} game
and think of it as ``racing against time''. For each $v \in V$, let $f_v(t)$ be the
probability generating function of the random variable ``number of turns'' until the end, 
if your current location is $v$.

In other words $f_v(t)$ is the formal power series whose coefficient (in its Maclaurin expansion) of $t^k$
is the probability of reaching $e$ from $v$ in exactly $k$ rounds.

The most important one is $f_1(t)$ where $1$ is the initial state, but it is is also
useful, during the game, if you are currently at vertex $v$, to know the probability
distribution (and hence the expectation) for the duration until the end of game.
Note that if our directed graph has cycles (like in the Chutes and Ladders game), the
$f_v(t)$ are {\bf not} polynomials but, as we will soon see, are {\bf rational functions} of $t$.

Suppose that you are currently at $v$. Then your next location is $u$ for some $u \in N(v)$, and your probability
of moving there in that round is $p_{vu}$. But by going there you have spent one round.
This leads to the {\bf linear equation}
$$
f_v(t) \, = \, t  \sum_{ u \in N(v)} \, p_{vu} f_u(t) \quad, \quad v \in V \backslash \{e\} \quad.
$$
$$
f_e(t)=1 \quad.
$$
The equation $f_e(t)=1$ follows from the fact that if you are at the absorbing state, the probability of getting there
in $0$ steps is $1$.

This gives us a system of $|V|$ linear equations with $|V|$ unknowns, that our computer can easily solve, for each
specific game. 

To see the list of length $99$ whose $i$-th entry is the probability generating function for the duration
of the ``Chutes and Ladders'' game (produced by {\it Winning Moves Games}, described above) see

{\tt http://sites.math.rutgers.edu/\~{}zeilberg/tokhniot/oSnakesAndLadders1.txt} \quad .

To see the list of length $99$ whose $i$-th entry is the probability generating function for the duration
of the ``Snakes and Ladders'' game (produced by {\it Cardinal Industries}, described above) see

{\tt http://sites.math.rutgers.edu/\~{}zeilberg/tokhniot/oSnakesAndLadders1a.txt} \quad .

By Cramer's rule it follows that the probability generating functions for durations $\{f_v(t)\}$ are
all {\bf rational functions}, and this leads to our first ``conceptual'' theorem.

{\bf Theorem 1}: For any finite game of pure chance, given by a Markov process as above, the
probability generating function for the random variable, ``number of rounds until the end''
starting at each particular location (in particular at the initial position) are
all rational functions with the {\bf same} denominator. Equivalently, the sequence
of probabilities are all $C$-recursive sequences satisfying the {\bf same} linear
recurrence equation with {\bf constant coefficients} (but of course with different initial values).
Furthermore, if the transition probabilities are rational (in particular, if the spinner is fair), then 
the coefficients of the numerators and denominators are {\bf  rational numbers}.

Once you have the probability generating function $f_v(t)$, you immediately get the
important quantity called the {\bf expectation} (average), the {\bf variance} and any desired moments.
The expectation is $t\frac{d}{dt} f_v(t)|_{t=1}$, and the $k$-th moment is
$(t\frac{d}{dt})^k f_v(t)|_{t=1}$. Since the derivative of a rational function is
rational, and the coefficients of the numerators and denominators are still rational numbers, we get.

{\bf Corollary 1.1}:  
For any finite game of pure chance, given by a Markov process as above, where
the transition probabilities are rational numbers, the expectation, variance, and higher moments
are certain specific {\bf rational numbers}.

The {\bf exact} value of the expected duration (from the starting location) for the {\it Winning Moves Games} version is
the {\bf rational number}

$$
{\frac {887878294805352403696983059454536608342612464186714311208985}{24341611817328043604468543532599921796578664874808501948232}}
\quad ,
$$

whose floating-point representation is  $36.475739629259028643943017\dots$ . So the solitaire game is expected to
last that long. The variance is also a rational number, and its square-root, the {\bf standard deviation} is
 $23.3564795406961083914719\dots$, note that it is fairly large. For higher moments, and more details, see the output file

{\tt http://sites.math.rutgers.edu/\~{}zeilberg/tokhniot/oSnakesAndLadders4.txt} \quad .

The {\bf exact} value of the expected duration (from the starting location) for the {\it Cardinal Industries} version is
the {\bf rational number}

$$
{\frac {2187389648884112026248013918019557612757259881117230420982560723}{74261533674379113296972151137271925497892767961930598965112832}}
\quad ,
$$

whose floating-point representation is  $29.45521780462205595995\dots$ . So the solitaire game is expected to
last that long. The variance is also a rational number, and its square-root, the {\bf standard deviation} is
$16.119642155096650734294068$, note that it is fairly large. For higher moments, and more details, see the output file

{\tt http://sites.math.rutgers.edu/\~{}zeilberg/tokhniot/oSnakesAndLadders5.txt} \quad .

{\bf Probability of Winning}

Suppose that there are  two players, where both players are currently at location $v$. In particular, if they are at the starting position.
Writing
$$
f_v(t)= \sum_{k=0}^{\infty} a_k t^k \quad ,
$$
the probability of the first player winning is
$$
\sum_{ 1 \leq k_2 \leq k_1 \, < \, \infty } a_{k_1} a_{k_2} \quad ,
$$
that is easily seen to equal 
$$
\frac{1}{2} \left ( 1+ \sum_{k=0}^{\infty} a_k^2 \right ) \quad .
$$

Since $\sum_{k=0}^{\infty} a_k t^k$ is a rational function in the variable $t$, it follows from elementary linear algebra (see [Z1]) that
 $\sum_{k=0}^{\infty} a_k^2 t^k$ is also a rational function, and if all the transition probabilities of our Markov process
describing the game  are rational numbers, the numerators and denominators are polynomials whose
coefficients are rational numbers, hence we get the surprising fact that the infinite convergent sum
$\sum_{k=0}^{\infty} a_k^2$ is a specific {\it rational number}.

We have demonstrated the ``two-player'', ``same starting position'' of the following theorem.
The general case can be proved similarly, and is left to the reader.

{\bf Theorem 2}:
For any finite game of pure chance, given by a Markov process as above, where
the transition probabilities are rational numbers, and there are $k$ players, 
currently all at the same location (in particular if they are at the initial position), or possibly different locations,
where they take turns moving, and the player to first reach the absorbing state is the winner, the 
winning probability of each of the players is {\bf always} a {\it specific} rational number,
that can be explicitly computed.

We were too lazy to find the actual rational numbers describing the
probability of the first player winning, at the very beginning of the game
(it is a long computation, since the probability generating function is complicated), but their
floating point approximations, for the above two versions of ``Snakes and Ladders'') discussed above,
are $0.5087744562$ and $0.5112590928$ respectively. See the above output files.

{\bf A short user's manual for the Maple package  {\tt SnakesAndLadders.txt} }

The material in the present chapter is implemented in the Maple package   {\tt SnakesAndLadders.txt},
available from the front of this article, or directly from

{\tt http://www.math.rutgers.edu/\~{}zeilberg/tokhniot/SnakesAndLadders.txt} \quad .

The main procedures are :

$\bullet$ {\tt GFD(P,t)}, that inputs a Markov process $P$ and a variable $t$ and outputs a list
of rational functions in $t$ corresponding to the functions $f_v(t)$ described above. In particular,
the first entry is the probability generating function for the duration of a solitaire game
starting at the beginning. First let us get a simple Markov process with the aid of the command
{\tt TMdieG}:

{\tt M:=TMdieG([[1,1/2],[2,1/2]],7,[[[1,3]],[[4,2]]] );}

you would give us the discrete Markov process

{\tt M=[[[2, 1/2], [3, 1/2]], [[2, 1/2], [3, 1/2]], [[2, 1/2], [5, 1/2]], [[5, 1/2], [6, 1/2]], [[6, 1/2], [7, 1/2]], [[7, 1]]]} \quad .

(This means that there are $6$ states, $1,2,3,4,5,6$, and the absorbing state is $7$. 
The probability of going from $1$ to $2$, and from $1$ to $3$ are both $\frac{1}{2}$,  The probability of going from
$2$ to $3$ is $\frac{1}{2}$, and from $2$ back to $2$ is also $\frac{1}{2}$ etc.).

Having gotten {\tt M}, typing  {\tt R:=GFD(M,t)[1];} would give the probability generating function, {\tt R}, of the duration starting at the first location, $1$.
It turns out to be
$$
\frac {{t}^{3} \left( 1+t \right) }{2(4-2\,t-{t}^{2})}
$$

$\bullet$ {\tt ProbAhead(R,t)}, inputs a rational function {\tt R}  that is the probability generating function for the game
(obtained thanks to {\tt GFD(M,t)}), and outputs the probability of the first player winning in the two-player version.
For example, calling

{\tt ProbAhead(R,t);}

would give the nice rational number $\frac{11}{20}$. Alas, if $R$ is very complicated, it may take a long time, so the approximate version

{\tt ProbAheadAppx(R,R, t,K);}

gives an approximation using the first {\tt K} terms, should be used for a sufficiently large {\tt K}.

For more details, explore the on-line Help (invoked by {\tt ezra();}).

{\bf Chapter 2: ``Infinite Families'' of ``Snakes and Ladders'' games, but with neither Snakes nor Ladders}.

In the previous chapter, we studied one game at a time. We now study  ``infinite'' families of games of pure chance
of the following kind.

The input is an {\it arbitrary} die (with an arbitrary, but finite, number of faces) each face with a certain positive number of dots,
and the die can be as loaded as one wishes. In other words, the input is an arbitrary probability distribution, let's call it $\P$,  on a finite
set of positive integers. If the die (or spinner) has $k$ faces with number of dots $i_1, \dots, i_k$, whose
respective probabilities are $p_1, \dots, p_k$ (of course $p_1+ \dots + p_k=1$). We denote it by
$$
\P = \{ \, [i_1, p_1] \, , \, [i_2, p_2] \, , \,  \dots \, , \,  [i_k, p_k] \, \} \quad .
$$
For example, for the familiar fair cubic die  the distribution is
$$
\P= [ \, [1,\frac{1}{6}] \, , \, [2,\frac{1}{6}] \, , \, [3,\frac{1}{6}] \, , \, [4,\frac{1}{6}] \, , \, [5,\frac{1}{6}] \, , \, [6,\frac{1}{6}] ] \quad .
$$
If you toss a loaded coin whose probability of Heads is $\frac{2}{3}$
and get one dollar if it lands on Heads, and $2$ dollars if it lands on Tails, the probability distribution is

$$
\P = [ \, [ \, 1 \, , \, \frac{2}{3} \, ] \, , \, [ \, 2 \, , \,  \frac{1}{3} \, ] \,] \quad .
$$

The other input is a positive integer $n$. In the Solitaire version you keep rolling the die, accumulating capital,
and end the game as soon as you have reached your goal of $\geq n$ dollars. In the Game version the players take turns
and the first person to reach the goal of $\geq n$ dollars wins the game.

Now for each {\it specific} positive integer $n$, this is a game similar to the one dealt with in Chapter $1$, only simpler,
since there are no cycles, so the probability generating functions are always polynomials, rather than rational functions.

But we want, having fixed the probability distribution $\P$, 
to get nice {\bf closed form expressions},  in terms of the {\bf symbol} $n$,
for the expectation, variance, and any desired higher moment for the  duration.

A rough estimate for the expected number of moves is $n/E[\P]$, where $E[\P]=\sum_{r=1}^{k} i_r p_r$ is the expected gain in one move, but
one can do much better as follows.

We need the {\bf grand generating} function, in $x$, say
$$
\sum_{n=0}^{\infty} F_n(t) x^n \quad,
$$
where $F_n(t)$ is the probability generating function for the $\P$-game with $n$ as goal. We clearly have
$$
\sum_{n=0}^{\infty} F_n(t) x^n \, = \,
\frac{\sum_{r=1}^{k} t\,(1+ x+ \dots + x^{i_r -1}) p_{i_r}}
{1-t(\sum_{r=1}^{k} p_r x^{i_r} ) } \quad .
$$

We would like to have explicit expressions for the expectation, variance, and higher moments {\bf in terms of $n$}.
Differentiating with respect to $t$, and plugging-in $t=1$ gives something of the form
$$
\sum_{n=0}^{\infty} F_n'(1) x^n \, = \, \frac{P(x)}{(1-x)^2 Q_1(x)} \quad,
$$
for {\it some} polynomial $P(x)$, and {\it some} polynomial $Q_1(x)$ whose roots are all larger than $1$ in absolute value. More generally
$$
\sum_{n=0}^{\infty} [(t \frac{d}{dt})^k) F_n(t)] \vert_{t=1} x^n \, = \, \frac{P(x)}{(1-x)^{k+1} Q_k(x)} \quad,
$$
for some polynomial $P(x)$ and {\it some} polynomial $Q_k(x)$ whose roots are all larger than $1$ in absolute value. 
Performing a partial fraction decomposition over the complex numbers leads to something of the form
$$
\sum_{n=0}^{\infty} [(t \frac{d}{dt})^k) F_n(t)] \vert_{t=1} x^n \, = \, \
\frac{A_0}{1-x} + \frac{A_2}{(1-x)^2} + \dots + \frac{A_k}{(1-x)^{k+1}} + \sum_{s=1}^m \frac{B_s(x)}{(x-\alpha_s)^{k+1}} \quad .
$$

Recalling that the coefficient of $x^n$ in $\frac{1}{(1-x)^{r+1}}$ is  ${{n+r} \choose {r}}$,
which is a polynomial in $n$ of degree $r$, and the coefficient of $x^n$ in $\frac{1}{(x-a)^{r+1}}$, with $|a|>1$ is
$o(1)$, we get our next `conceptual' theorem.

{\bf Theorem 3}: Given an arbitrary finite probability distribution $\P$ as above, 
except for exponentially small terms, the {\bf average}  of the random variable, ``number of rounds it takes to
reach $n$'', is a  {\bf polynomial} in $n$ of degree $1$. Furthermore, for the higher moments, the $k$-th moment of that random variable
(and hence the $k$-th moment about the mean) is  a {\bf polynomial} in $n$ of degree $\leq k$, that can
be explicitly computed.

Here are some sample results.

{\bf Proposition 1}: For $\P \,=\, [\,[1, \frac{1}{2}] \,, \, [2, \frac{1}{2}]]$, up to exponentially small contributions,
for the random variable, 

`duration until you get $\geq n$ for the first time', we have

$\bullet$ The expectation is $\frac{2}{3}\,n+ \frac{2}{9} + o(1)$ \quad .

$\bullet$ The variance  is $\frac{2}{27}\,n+ \frac{2}{81} + o(1)$ \quad .

$\bullet$ The third moment about the mean   is $\frac{2}{81}\,n -  \frac{26}{729} + o(1)$ \quad .

$\bullet$ The fourth moment about the mean   is  ${\frac {4}{243}}\,{n}^{2}  \,+\,{\frac {2}{243}}\,n \, - \, {\frac {62}{2187}}+ o(1)$ \quad .

For the $5$-th through the $10$-th moments, see the output file

{\tt http://sites.math.rutgers.edu/\~{}zeilberg/tokhniot/oPosPileGames1.txt} \quad .

More generally, for a loaded coin, where the probability of $1$ is $p$ and the probability of $2$ is $1-p$, we have the
next proposition.

{\bf Proposition 2}: For $\P \,= \, [ \, [1, p] \, , \, [2, 1-p] \, ]$, up to exponentially small contributions,
for the random variable, 

`duration until you get $\geq n$ for the first time', we have

$\bullet$ The expectation is 
$$
{\frac {1}{2-p}}\, \cdot n \, + \, {\frac {1-p}{ \left( 2- p\right) ^{2}}} \, + \, o(1) \quad .
$$

$\bullet$ The variance  is 
$$
\frac{p(1-p)}{(2-p)^3} \cdot n \, - \,\frac{(1-p)(p^2+p-1)}{(2-p)^4} \, + \, o(1)
$$

$\bullet$ The third moment about the mean   is 
$$
-{\frac {p \left( 1-p \right)  \left( {p}^{2}+2\,p-2 \right)}{ \left( 2-p \right) ^{5}}} \cdot n
\, + \, {\frac {{p}^{2} \left( 1-p \right)  \left( {p}^{2}+7\,p-7 \right) }{ \left( 2-p \right) ^{6}}}
\, + \, o(1) \quad .
$$
$\bullet$ The fourth moment about the mean   is  

$$
\frac{3\,p^2\,(1-p)^2}{(2-p)^6} \cdot n^2
\, + \, {\frac {p \left( 1-p \right)  \left( {p}^{4}+16\,{p}^{3}-6\,{p}^{2}-20\,p+10 \right) }{ \left( 2-p \right) ^{7}}} \cdot n
$$
$$
+ \, {\frac { \left( -1+p \right)  \left( {p}^{6}+26\,{p}^{5}+12\,{p}^{4}-75\,{p}^{3}+35\,{p}^{2}+3\,p-1 \right) }{ \left( -2+p \right) ^{8}}}
 \, + \, o(1) \quad .
$$

For the $5$-th through the $10$-th moments, see the output file

{\tt http://sites.math.rutgers.edu/\~{}zeilberg/tokhniot/oPosPileGames2.txt} \quad .

Maple can easily generate such proposition for each {\it specific} die, but what about an `infinite' family of dice?
For the `infinite' family of $k$-faced fair dice, for each case, 
$k=2,3,4, \dots$, the algorithm that we used can crank out explicit expressions, in $n$,
(up to exponentially small terms) for the expectation, variance, and any desired finite moment, but it can't do it
(at least not with the present method) for {\it symbolic} $k$, i.e. all $k$ at once.
But it is possible to show that these quantities for the $k$-faced fair die , i.e. for
$$
\P \, = \,[\, [1, \frac{1}{k}] \,, \, [2 , \frac{1}{k}] \,, \, \dots \, , \,  [k , \frac{1}{k}] \, ] \quad,
$$
in addition to being polynomials in $n$, are also {\bf rational functions} in $k$. Being experimental mathematicians, we collected
enough data for several $k$, and then ``fitted'' it with a rational function. leading to the next impressive, {\bf computer-generated}
proposition.

{\bf Proposition 3}: For {\it any} positive integer, and any fair $k$-sided die,
for the random variable,  `duration until you get $\geq n$ for the first time', we have, up to exponentially small terms  

$\bullet$ The expectation is 
$$
{\frac {2}{k+1}} \cdot n \, + \,\frac {2(k-1)}{3(k+1)}   \, + \, o(1)\quad .
$$

$\bullet$ The variance  is 
$$
\frac{2(k-1)}{3(k+1)^2} \cdot n  \, + \, \frac{2(k-1)^2}{9(k+1)^2} \,+   \, o(1) \quad .
$$

$\bullet$ The third moment about the mean is
$$
\frac{2}{3} \,{\frac { \left( k-1 \right) ^{2}}{ \left( k+1 \right) ^{3}}} \cdot n 
\,+\,{\frac {2}{135}}\,{\frac { \left( k-1 \right)  \left( k-7 \right)  \left( 7\,k-1 \right) }{ \left( k+1 \right) ^{3}}}
\, + \, o(1) \quad .
$$

$\bullet$ The fourth moment about the mean is
$$
\frac{4}{3}\,{\frac { \left( k-1 \right) ^{2}}{ \left( k+1 \right) ^{4}}} \cdot n^2 \, 
\, + \, \frac{2}{15} \,{\frac { \left( k-1 \right)  \left( 13\,{k}^{2}-30\,k+13 \right)}{ \left( k+1 \right) ^{4}}} \cdot n
\, + \, {\frac {2}{135}}\,{\frac { \left( 13\,{k}^{2}-110\,k+13 \right)  \left( k-1 \right) ^{2}}{ \left( k+1 \right) ^{4}}}
+ \, o(1) \quad 
$$

For the fifth and sixth moments, see the output file

{\tt http://sites.math.rutgers.edu/\~{}zeilberg/tokhniot/oPosPileGames4.txt} \quad .

\vfill
\eject

{\bf Probability of Winning}

So far we considered the solitaire game, but now let's turn it into a two-player game, where the players
take turns, and whoever reaches the goal $n$ first, is the winner. For each {\it specific} $n$,
we know from the previous chapter, that it is a specific rational number (provided the probabilities in
$\P$ are rational), but what can you say about the {\it sequence}, let's call it $f(n)$, of the first
player winning?

By {\bf Wilf-Zeilberger algorithmic proof theory} [PWZ][Z3], the {\it double} sequence,
let's call it $b_{k,n}$, the coefficient of $t^k x^n$ in the grand-generating function above,
is {\it holonomic} in both $n$ and $k$, i.e. satisfies linear recurrence equations with polynomial coefficients
in both the $n$ and $k$ variables. It also follows from that theory that the sum of the squares
$\sum_{k=1}^{\infty} b_{k,n}^2$, let's call it $a(n)$, is holonomic ($P$-recursive)
in the surviving variable $n$, i.e.
$a(n)$ satisfies some specific linear recurrence equation with polynomial coefficients, that enables
a very fast computation of many terms, once that recurrence is known. 

This brings us to the next `conceptual' theorem.

{\bf Theorem 4}: Given an arbitrary finite probability distribution $\P$ as above, in the two player-game
where players take turns and the first to reach $n$ is the winner, let $f(n)$ be the probability
that the first player wins. Then $f(n)=(1+a(n))/2$, where
the sequence $a(n)$ (and hence, also $f(n)$) is $P$-recursive, i.e. satisfies a linear
recurrence equation with polynomial coefficients.

While there exist algorithms to do this {\it ab initio}, it is much more efficient to crank-out enough terms
of the desired sequence and use {\it undetermined coefficients} to discover the recurrence, in the
spirit of {\it experimental mathematics}.

This brings us to the next computer-generated proposition.

{\bf Proposition 4}: if two players take turns tossing a fair coin and get one dollar if it is Heads and two dollars
if it is Tails, and the first to reach $n$ dollars is the winner, the probability of the  player who goes first to
win the game is $\frac{1}{2} (1+ a(n))$, where $a(n)$ satisfies the linear recurrence
$$
a \left( n \right) \,= \, \frac{1}{2}\,{\frac { \left( 3\,n-1 \right)  \left( n-3 \right) }{n \left( 3\,n-7 \right) }} \cdot  a(n-1)
+ \frac{1}{16}\,{\frac { \left( 21\,{n}^
{2}-67\,n+62 \right) }{n \left( 3\,n-7 \right) }} \cdot  a(n-2)
$$
$$
\,+\,  \frac{1}{16} \,{\frac { \left( 6\,{n}^{2}-17\,n+2 \right) }{n \left( 3\,n
-7 \right) }}  \cdot a( n-3) \,- \, \frac{1}{16} \,{\frac { \left( n-4 \right)  \left( 3\,n-4 \right) }{n \left( 3\,n-7 \right) }} 
\cdot  a( n-4)\quad ,
$$
subject to the initial conditions
$$
a(1) = 1 \,, \,  a(2) = \frac{1}{2} \, , \, a(3) = \frac{5}{8} \, , \,  a(4) = \frac{15}{32} \quad .
$$

Using this recurrence  it follows that the probability of the first player
to reach $n=1000$ first is $(1+a(1000))/2 \, = \, 0.516384982\dots$.

{\bf Comment}: While for the general case it is much easier to use the `guessing' way (that is easily made rigorous by invoking
general theorems), in this simple case, where the probability of ending after exactly $k$ rounds, if the goal is $n$,
is easily seen to be given by the closed-form expression
$$
b_{k,n} \, = \, {\frac {{k-1\choose n-k} \left( 3\,k-n \right) }{ \left( 2\,k-n \right) {2}^{k}}} \quad ,
$$
to get a recurrence satisfied by the sum of squares, one can use the celebrated Zeilberger algorithm (see [PWZ]), implemented
in Maple. Just type:

{\tt
 ope:=SumTools[Hypergeometric][Zeilberger]
} \hfill\break
{\tt
(( binomial(k-1,n-k)*(3*k-n)/(2*k-n)/2**k)**2,n,k,N)[1];
}

Followed by (to make it look nicer)

{\tt
add(factor(coeff(ope,N,i)/coeff(ope,N,4))*N**i,i=0..4); \quad ,
}

that is equivalent to Proposition 4. (Recall that $N$ is the forward shift operator in $n$: $Na(n):=a(n+1)$).

For the more general case where the probability of winning a dollar is $p$, rather than $\frac{1}{2}$, see the output tile

{\tt http://sites.math.rutgers.edu/\~{}zeilberg/tokhniot/oPosPileGames6.txt } \quad .

The front of this article contains a few other such propositions, and readers can create their own.

{\bf Chapter 3: Games of Pure Chance  Generated by Gambler's Ruin with Unlimited Credit: The Fuss-Catalan case}

In Chapter 2, we only allowed  positive steps. Now
we will also allow negative steps, and treat  games that may be viewed as a ``gambler's ruin with infinite credit"
with an {\it arbitrary} `die'. In the next chapter we will treat the case of a {\it general}  die,
while in this chapter we only consider two-faced dice, where one of the faces is marked $1$ and the other marked
$-k$, and the more difficult case where one of the faces is marked $-1$ and the other marked $k$. We will
start with their intersection, the very classical case of $\{-1,1\}$, treated in Feller's classic [F].
However, even in this case we will be able to go beyond Feller, since he did not use a computer.

The general set-up, to be considered in full generality in Chapter 4 is as follows.

On the discrete line, you start at the origin $x=0$, and there is a fixed allowed set of steps consisting of both positive and negative integers
and a probability distribution on them, let's call it $\P$. 
You are allowed to go as far left as possible (i.e. you can owe as much as necessary).
At each round, you roll the $\P$ die, and move accordingly. You win as soon as you reach a location $\geq1$, or more generally
when you reach a location $\geq n$. In other words, your goal is to exit the casino with at least one dollar (or more generally, at least $n$
dollars). In the two-player (or multi-player) version, the players take turns rolling the $\P$ die, and whoever
achieves the goal first is declared the winner. As before, we will first discuss the solitaire game, where the
goal is to reach it as soon as possible.

{\bf Classical Gambling:  Winning a dollar or losing a dollar}

Let's start with the simplest, most classical case, of simple random walk, where you start with $0$ dollars,
and at each round you win a dollar with probability $p$ and lose a dollar with probability $1-p$.
The expected gain at each individual round is $p\cdot 1+(1-p) \cdot (-1)=2p-1$, so if $p>\frac{1}{2}$, then
sooner or later you will reach your goal of owning $\geq 1$ dollars. If $p<\frac{1}{2}$, then you may never make it, sliding down to
infinite debt. In the border-line case of a {\bf fair} coin, $p=\frac{1}{2}$, as we will soon see, you are
also guaranteed to `eventually' be in possession of $1$ dollar (and more generally, $n$ dollars for each $n>0$,
as big as you wish). Alas, as we will also soon see, the expected time until that happens is infinite, and since
life is finite, there is a good chance that when you will pass away, your heirs will have a huge debt.

{\bf Analyzing Gambling  histories}

For typographical clarity, let's denote $-1$ by $\1$.

Our {\bf alphabet} is $\{-1,1\}$=$\{\1,1\}$. A `gambling history' consists of a {\bf word} that
ends in $1$, whose sum is $1$, and whose proper partial sums are all non-positive.
Obviously the length of such a game is odd.

If you are really lucky, you exit after one step, since you won a dollar right away.

If you lost a dollar at the first round, you can recover at the second round, and then  win a dollar at the third round. Etc.

For the sake of clarity and concreteness, let's list the first few `histories'.

Length $1$: $\{\, 1 \, \}$. Probability $=p$.

Length $3$: $\{\, \1 \,1 \, 1 \, \}$. Probability $=p^2\,(1-p)$.

Length $5$: $\{\1\, 1 \, \1 \, 1 \, 1  \quad , \quad \1 \,\1 \,1 \,1 \,1 \, \}$. Probability $2 \cdot p^3\,(1-p)^2$.

Length $7$: 
$$
\{
\1 \, 1 \, \1 \, 1 \, \1 \, 1 \, 1 \quad, \quad
\1 \, 1 \,   \1 \, \1 \, 1 \,  1 \, 1 \, \quad, \quad
\1 \, \1 \, 1 \, 1 \, \1 \, 1 \, 1 \, \quad, \quad
\1 \, \1 \, 1 \, \1 \, 1 \, 1 \, 1  \quad , \quad 
\1 \, \1 \, \1 \, 1 \, 1 \, 1 \, 1 \,
\} \quad ,
$$
with probability $5 \cdot p^3\,(1-p)^2$.

 It is useful, for humans, to visualize such a history as
a lattice path in the discrete plane starting at $(0,0)$ where $\1$ corresponds to a step $(1,-1)$ and
$1$ corresponds to a step $(1,1)$. For example, the word (gambling history)
$$
\1 \, \1 \, 1 \, 1 \, \1 \, 1 \, 1 \quad,
$$
corresponds to the walk
$$
(0,0) \rightarrow (1,-1)
\rightarrow (2,-2)
\rightarrow (3,-1)
\rightarrow (4,0)
\rightarrow (5,-1)
\rightarrow (6,0)
\rightarrow (7,1) \quad .
$$

Let's study the {\it anatomy} of such histories, or equivalently, paths . Obviously they are all of odd length, and they all end with $1$.
So we can write, for {\it any} history $W$
$$
W=U\, 1 \quad,
$$
where $U$ is a word that sums to $0$, all whose partial sums are non-positive. Such words are called {\it Dyck words}.

Let's analyze such a Dyck word $U$ or rather its corresponding path from $(0,0)$ to $(2n,0)$, say.
Of course, it may be the empty word, but if it is not,
let $(2r,0)$ $0<r\leq n$ be the {\bf first} time that it hits the $x$-axis. Then we can write
$$
U= U_1 \, U_2 \quad ,
$$
where $U_2$ is another word of that kind (of length $2n-2r$), but $U_1$, consisting of the first $2r$ letters of $U$, has the special property that all its partial sums
(except the $0$-th and the last) are {\bf strictly negative}, or in terms of its path, except for its
starting and ending points, they lie {\bf strictly} below the $x$-axis. Such a word must {\bf necessarily}
start with a $\1$ and end with a $1$, and may be written as $\1 U_3 1$, where $U_3$ is an arbitrary
Dyck word.  Conversely, for any Dyck word $U_3$, $\1 \, U_3 \, 1$ corresponds to such a `strictly below the $x$-axis' path.
So we have the (context-free) {\it grammar}
$$
U= EmptyWord \, \vee \, \1\, U \, 1 U \quad,
\eqno(DyckGrammar)
$$
where now $U$ stands for `an arbitrary Dyck word'.

let $z_1$ and $z_{-1}$ be {\bf commuting} variables.

For any word $u=u_1 \dots u_m$, let the {\it weight} of $u$ be $z_{u_1} \cdots z_{u_m}$. For example, 

$$
weight(\1 \, \1 \, \1 \, 1 \, 1 \, \1 \, 1 \, 1)=
z_{-1}z_{-1}z_{-1}z_1 z_1 z_1 z_{-1} z_1 z_1
\, = \, z_{-1}^4 z_1^5 \quad .
$$

Let $F(z_{-1},z_1)$ be the {\bf weight enumerator} of the set of Dyck words,
i.e. the {\it sum} of all the weights of all these words, a certain {\bf formal power series} in $z_{-1},z_1$.

Obviously the weight of the empty word is $1$ (the empty product), hence applying $weight$ to $(DyckGrammar)$, we get
the {\bf quadratic equation}
$$
F = 1 \,+ \, z_{-1}\,F \, z_1 \, F \quad .
$$

Abbreviating $X=z_{-1}\,z_1$, we get
$$
F=1 \,+ \, X \, F^2 \quad .
$$
Recalling what we learned in seventh grade (or what the Babylonians knew more than $3000$ years ago), we
can express $F$ {\bf explicitly}
$$
F= \frac{1 -\sqrt{1-4X}}{2X} \quad .
$$
Recalling what we learned in $12$-th grade (or what Isaac Newton knew more that $300$ years ago) we can write
$$
F= \sum_{m=0}^{\infty} \frac{(2\,m)!}{m!\,(m+1)!} X^m \quad ,
$$
implying the fact that the number of  Dyck paths  of length $2m$
is the super-famous {\bf Catalan} number $C_m \,= \, \frac{(2m)!}{m!\,(m+1)!}$, that is the subject of
Richard Stanley's modern classic [St], and the most popular sequence, {\tt A108}, in the great OEIS [Sl].

{\eightrm
The above is the standard, very boring proof of that famous fact. We know at least a dozen proofs, some of them are given in [St].
Here is one of our favorite proofs due to Aryeh Dvoretzky and Theodore Motzkin [DM].

The fact that the number of Dyck paths of length $2m$ equals the Catalan number $C_m$
is equivalent the fact that the number of words in $\{1,-1\}$ of length $2m+1$ whose sum is $1$ and all
whose proper-partial sums are non-positive is $C_m$. 
Every word of length $2m+1$ in $\{-1,1\}$ that adds up to $1$ has $m+1$ `$1$' and $m$ `$\1$'. There
are ${{2\,m \, + \, 1} \choose {m}}$ such words. The $2\,m \, + \, 1$ {\it cyclic shifts} of each such word are all {\bf different} (why?), and
exactly one of them has the property that its partial sums are all non-positive (why?). Hence the number
of gambling histories that we are interested in is $\frac{1}{2m\,+\,1} \cdot {{2m+1} \choose {m}}=C_m$.
}

{\bf Enter Probability}

So far what we did was {\it enumerative combinatorics}. 
We found out that the weight-enumerator of the set of Dyck words is
$$
\frac{1 -\sqrt{1-4z_{-1}z_1}}{2\,z_{-1}z_1} \quad,
$$
and hence the weight enumerator of words in $\{-1,1\}$ that add-up to $1$,  and such that all their proper partial sums are $\leq 0$, is
$z_1$ times that, i.e.
$$
\frac{1 -\sqrt{1-4z_{-1}z_1}}{2\,z_{-1}} \quad.
$$
Assume that each round in the gambling game is {\bf independent} of the other ones, and for each of them the probability of
winning a dollar is $p$, and hence of losing a dollar is $1-p$. Plugging-in $z_{-1}=(1-p)\,t$, $z_1=p\,t$, in the above
explicit enumerating generating function, we get the following human-generated, well-known (see [F])  proposition.

{\bf Proposition 5}: The {\bf probability generating function} of the random variable `numer of rounds it takes until
the first time you have one dollar', if you start with $0$ dollars and at each round you win a dollar
with probability $p$ and lose a dollar with probability $1-p$, let's call it $g(t)$, is
$$
g(t)=\frac{1 -\sqrt{1-4\,(1-p)\,p t^2}}{2\,(1-p)\,t} \quad.
$$

So far all our power series were {\it formal}, but it is easy to see that if $p \geq \frac{1}{2}$ then plugging-in $t=1$
leads to a convergent series, that sums-up to $1$, in agreement with the obvious fact that if $p>\frac{1}{2}$ sooner
or later you will succeed, and the slightly less obvious fact that it is still true when $p=\frac{1}{2}$. If
$p < \frac{1}{2}$, then we must take the other sign of the square-root, leading to the classical and well-known fact that
the probability of one day having one dollar in your possession is $\frac{p}{1-p}$.

More generally, suppose that your goal in life is not just to exit the casino with one dollar, but you want to
make $n$ dollars. Since each additional dollar is yet another 1-dollar game, we immediately get.

{\bf Proposition 5'}: The {\bf probability generating function} 
of the random variable `numer of rounds it takes until
the first time you have $n$ dollars', if you start with $0$ dollars and at each round you win a dollar
with probability $p$ and lose a dollar with probability $1-p$,  is given by
$$
\left ( \frac{1 -\sqrt{1-4\,(1-p)\,p t^2}}{2\,(1-p)\,t} \right )^n\quad.
$$

From now let's assume that $p \geq \frac{1}{2}$. To get the {\bf expected duration} we can sill do it by hand,
find $(g(t)^n)' \,= \, ng(t)^{n-1} g'(t)$, then compute $g'(t)$, plug-in $t=1$ and simplify, getting
that the expectation is $\frac{n}{2p-1}$.

For the $k$-th  moment, we compute $(t \frac{d}{dt})^k (g(t)^n)$, plug-in $t=1$, and simplify, expressing all higher derivatives
of $g(t)$ in terms of $g(t)$ and $t$, followed by substituting $t=1$.

An even better way, that would be the only way later on when we do the
general gambling caes, is to use {\it implicit differentiation}, using the relation
$$
f(t) \,= \, 1+p\,(1-p)\,t^2 f(t)^2 \quad,
$$
and its implied relation for $g(t)=p\,t\,f(t)$.

It turns out that if you use the explicit expression 
$g(t)=\frac{1 -\sqrt{1-4\,(1-p)\,p t^2}}{2\,(1-p)\,t}$ all the radicals disappear, 
and if you use implicit differentiation, and then plug-in $t=1$,
you never have to divide $0$ by $0$, so either way you would get that all the moments are {\bf polynomials} in $n$ and rational functions in $p$.
In particular, if $p$ is a rational number, then they are all also rational numbers.
The expectation, is $\frac{n}{2p-1}$.

For higher moments, We get the following computer-generated proposition.

{\bf Proposition 6:} Let $X_{n,p}$ be the random variable ``Number of rounds until reaching $n$ dollars for the first time" in a gambling game
where the probability of winning a dollar is $p$ and of losing a dollar is $1-p$. Assume that $p > \frac{1}{2}$.
We have

$$
E[X_{n,p}] \, = \, \frac{n}{2p-1} \quad .
$$
$$
Var[X_{n,p}] \, = \, \frac{4\,n\,p\,(1-p)}{(2p-1)^3} \quad .
$$

The {\it skewness} (aka scaled third moment about the mean) is
$$
\alpha_3 [X_{n,p}] \, = \,
 \left( -2\,{p}^{2}+2\,p+1 \right)  \left( -1+2\,p \right) ^{-2}{\frac {1}{\sqrt {-{\frac {np \left( -1+p \right) }{ \left( -1+2\,p \right) ^{3}}}}}}
\quad .
$$

The {\it kurtosis} (aka scaled fourth moment about the mean) is
$$
\alpha_4 [X_{n,p}] \, = \,
{\frac {-4\,{p}^{4}+ \left( 6\,n+8 \right) {p}^{3}+ \left( -9\,n+6 \right) {p}^{2}+ \left( 3\,n-10 \right) p-1}{np \left( -1+p \right)  \left( -1+2\,p
 \right) }} \quad .
$$

For the $5$-th through $10$-th scaled moments, see the output file

{\tt http://sites.math.rutgers.edu/\~{}zeilberg/tokhniot/oGenPileGames1.txt} \quad .

{\bf The Two Player version for the $(1,-1)$ case}

Using Lagrange inversion (see [Z4] for a lucid statement and proof) or otherwise, it is easy to see that the
probability of reaching $m$ dollars for the first time after exactly $n$ rounds,
in a solitaire game where the
probability of winning a dollar is $p$ and the probability of losing a dollar is $1-p$, let's call it $b_{n,m}$ is
$$
b_{n,m} \, = \,
{\frac {m \left( 2\,n+m-1 \right) !\,{p}^{n+m} \left( 1-p \right) ^{n}}{n!\, \left( n+m \right) !}} \quad.
$$

Suppose that two players take turns and whoever reaches $m$ dollars first is declared the winner. As before, the
probability of winning the game for the player whose turn is to move is $a(m)=(1+f(m))/2$, where
$$
f(m)= \sum_{n=1}^{\infty} b_{n,m}^2 \quad .
$$
Using the Zeilberger algorithm once again we have the next computer-generated proposition.

{\bf Proposition 7}: In the two player version  game with a {\bf fair} coin, i.e. the probability of
winning a dollar and losing a dollar are both $\frac{1}{2}$, the winning probability of the player whose turn is to move 
is $(1+f(m))/2$ where $f(m)$ satisfies the second-order recurrence
$$
 \left( 2\,{m}^{2}+5\,m+2 \right) f \left( m+2 \right) + \left( -12\,{m}^{2}-24\,m-10 \right) f \left( m+1 \right) 
+ \left( 2\,{m}^{2}+3\,m \right) f \left( m \right) \, =\, -\frac{8}{\pi} \quad ,
$$
subject to the initial conditions
$$
f \left( 1 \right) =-{\frac {-4+\pi }{\pi }} \quad, \quad 
f \left( 2 \right) =-{\frac {-16+5\,\pi }{\pi }} \quad .
$$

For the loaded case, where $p>\frac{1}{2}$, we have the next proposition.

{\bf Proposition 8}: In the two player version  game with the probability of
winning a dollar is $p$ and losing a dollar is $1-p$ , provided $\frac{1}{2}<p<1$,
the winning probability of the player whose turn is to move
is $(1+f(m))/2$ where $f(m)$ satisfies the fourth-order recurrence
$$
m \left( -1+p \right) ^{4} \left( m-3 \right) f \left( m \right) 
- \left( -1+p \right) ^{2} \left( 2\,{m}^{2}-7\,m+4 \right) f \left( m-1 \right) 
$$
$$
+ \left( -2\,{m}^{2}{p}^{4}+4\,{m}^{2}{p}^{3}+8\,m{p}^{4}-2\,{m}^{2}{p}^{2}-16\,m{p}^{3}-4\,{p}^{4}+8\,m{p}^{2}+8\,{p}^{3}+{m}^{2}-4\,{p}^{2}-4\,m+4 \right) 
f \left( m-2 \right) 
$$
$$
-{p}^{2} \left( 2\,{m}^{2}-9\,m+8 \right) f \left( m-3 \right) +{p}^{4} \left( m-1 \right)  \left( m-4 \right) f \left( m-4 \right) =0
\quad ,
$$
subject to the appropriate initial conditions.

{\bf  Winning a dollar or losing k dollars}

Now let's generalize to the gambling game where, as before, you start with a capital of $0$  dollars,
but now at each round you {\bf win}  a dollar with probability $p$ or {\bf lose} $k$ dollars
with probability $1-p$, and the game ends as soon as you own $1$ dollar. Very soon we will treat the more general
case where the goal is to exit with  $m$ dollars, but for now let's consider the case of $m=1$.

In order to guarantee that the game ends, the expected gain of a single round, 
$p\cdot 1 - (1-p)\cdot k \, = \, (k+1)\,p \,- \, k$ should be positive.
So we will assume that $p > \frac{k}{k+1}$. In the border-line case $p=\frac{k}{k+1}$ the game still ends with
probability $1$, but its expected duration is infinite.

Now the {\bf alphabet} is $\{1, -k\}$, and we will try to adapt the above argument that worked for the classical case. 
Let's abbreviate $\k:=-k$. Now the steps are $(1,1)$ and $(1,-k)$.

Let's study the {\it anatomy} of such words (histories) or, equivalently,  paths. Obviously 
all these words are of length $n(k+1)+1$, for some non-negative integer $n$, and they all end with $1$.
So we can write, for {\it any} history $W$,
$$
W=U\, 1 \quad,
$$
where $U$ is a word that sums to $0$, all whose partial sums are non-positive. we will call such words  {\it $(1,-k)$-Dyck words}.

Let's analyze such a $(1,-k)$-Dyck word $U$ or rather its corresponding path from $(0,0)$ to $((k+1)n,0)$, say.
Of course, it may be the empty word, but if it is not,
let $(r(k+1),0)$ $0<r\leq n$ be the {\bf first} time that it hits the $x$-axis. Then we can write
$$
U= U_1 \, U_2 \quad ,
$$
where $U_2$ is another  arbitrary $(1,-k)$-Dyck word, but $U_1$ has the special property that all its partial sums
(except the $0$-th and the last) are {\bf strictly negative}, or in terms of its path, except for its
starting and ending points, they lie {\bf strictly} below the $x$-axis. Such a word must {\bf necessarily}
start with a $\k$ and end with a $1$, but to recover the `debt' of $k$, must regain these lost $k$ dollars, one
dollar at a time, so it
 may be written as $\k (U_3 1)^k$, where $U_3$ is an arbitrary
 $(1,-k)$-Dyck  word.  Conversely, for any such word $U_3$, $\k \, (U_3 \, 1)^k$ is such a strictly below the $x$-axis word.
So we have the (context-free) {\it grammar}
$$
U= EmptyWord \, \vee \, \k \, (U \, 1)^k U \quad,
\eqno((1,-k)-DyckGrammar)
$$
where now $U$ stands for `an arbitrary  $(1,-k)$-Dyck word'.

Let $F(z_{-k},z_1)$ be the {\it weight-enumerator} for all such words. Applying the {\it weight} operation,
we get that   $F=F(z_{-k},z_1)$ satisfies 
$$
F= 1 \,+\, (z_{-k}z_1^k) \, F^{k+1} \quad .
$$
Abbreviating $X := z_{-k} \, z_1^k$, this can be written
$$
F= 1\,+\, X \, F^{k+1} \quad .
$$
When $k=2$ and $k=3$, we can solve these equations `explicitly' using `radicals', thanks to Cardano and Ferrari, but
thanks to Abel, Ruffini, and Galois we know that we can {\bf not} do it for $k \geq 4$. Even the `explicit' solutions for $k=2$ and $k=3$ are
not very useful. On the other hand, thanks to {\it Lagrange inversion} (see,e.g. [Z4]) we can find the
Maclaurin expansion explicitly.
$$
F(X) \, = \, \sum_{m=0}^{\infty} \frac{((k+1)\,m)!}{m!\,(km \,+ \, 1)!} X^m \quad ,
$$
featuring the {\bf Fuss-Catalan} numbers   $C_{k,m}=\frac{((k+1)\,m)!}{m!\,(km\,+\,1)!}$ .

It follows that the weight-enumerator of words in $\{-k,1\}$ that add-up to $1$, and such that the proper-partial sums
are all non-positive is $F(z_{-k} z_1^k)\, z_1$, since the last letter must be $1$.

{\eightrm
Equivalently (and that's is our actual object of interest) the number of words with $m$ `$-k$' and $mk+1$ `$1$'
whose proper-partial sums are all non-positive equals the  Fuss-Catalan number $C_{k,m}$. This can
be also proved by adapting the  [DM] proof. There are ${{mk+1+m} \choose {m}}$ words altogether, and
for each of these its $mk+1+m$ cyclic shifts are all different, and exactly one of them is a `good' word,
hence there are  $\frac{1}{mk+1+m} \,{{mk+1+m} \choose {m}}=C_{k,m}$ such words.
}

Since, in order to exit with  $n$ dollars , we must gain one dollar, $n$  times, the weight-enumerator of
words that reach $n$ for the first time is $(F(z_{-k} z_1^k)\, z_1)^n$.

So far we did enumerative combinatorics. To convert it to probability, we plug-in the above
$z_1\,=\, p\,t$ and $z_{-k}=(1-p)\,t$. Using implicit differentiation, we can compute the
expectation, variance, and higher moments. Since in this case we do not encounter $0/0$, all the
moments are {\bf rational} functions of $p$. In particular, if the number $p$ is rational, all
the quantities are rational numbers.

Using implicit differentiation, for {\bf symbolic} $k$ and {\bf symbolic} $p$ and {\bf symbolic} $n$, 
our beloved computer generated the next proposition.

{\bf Proposition 9}: Suppose  that at each round, you win a dollar with probability $p$ and lose $k$
dollars with probability $1-p$, and you quit as soon as you reach $n$ dollars. If $p>k/(k+1)$, then, of course, sooner or 
later you will reach your goal. How long should it take?
Denote by $X_{n,k,p}$ the random variable, `number of moves until reaching $n$ dollars'. We have the following facts.

Let $g(t)$ be the formal power series, in $t$, satisfying the algebraic equation
$$
g \left( t \right) -1-{p}^{k} \left( 1-p \right) {t}^{k+1} \, g(t)^{k+1} \,=\, 0 \quad .
$$
The probability generating function of $X_{n,k,p}$ is
$$
(\,p \, t \, g(t) \,) ^n \quad .
$$
By implicit differentiation, followed by substituting $t=1$, we can compute any desired derivative, and hence the expectation, variance, and higher moments.
We have
$$
E[X_{n,k,p}] \, = \, \frac{n}{(p-1)k+p} \quad ,
$$
[as {\it expected} (npi), since the expected gain in one move is $(p-1)k+p$ ]. The variance is given by
$$
Var[X_{n,k,p}] 
\, = \, {\frac {np \left( k+1 \right) ^{2} \left( p-1 \right) }{ \left(  \left( 1-p \right) k+p \right) ^{3}}}
\quad .
$$
The {\it skewness} (aka `third scaled-moment about the mean') is
$$
\alpha_3 [X_{n,k,p}]  \, = \,
- \left( k+1 \right)  \left( k{p}^{2}+{p}^{2}-k-2\,p \right)  \left( kp-k+p \right) ^{-2}{\frac {1}{\sqrt {-{\frac {np \left( k+1 \right) ^{2} \left( p-1
 \right) }{ \left(  \left( p-1 \right) k+p \right) ^{3}}}}}} \quad .
$$
The {\it kurtosis} (aka `fourth scaled-moment about the mean') is
$$
\alpha_4 [X_{n,k,p}]  \, = \,
$$
$$
        {\frac {- \left( k+1 \right) ^{2}{p}^{4}-2\, \left( k+1 \right)  \left( k- \frac{3}{2}\,n-3 \right) {p}^{3}
+ \left( 6\,{k}^{2}+ \left( -6\,n+6 \right) k-3\,n-6
 \right) {p}^{2}-2\,k \left( k- \frac{3}{2} \,n+4 \right) p-{k}^{2}}{np \left( p-1 \right)  \left( p \left( k+1 \right) -k \right) }}
\quad .
$$

For the scaled fifth and sixth moments, see the output file

{\tt http://sites.math.rutgers.edu/~zeilberg/tokhniot/oGenPileGames2.txt} \quad .

{\bf The Two Player version for the $(1,-k)$ case}

Since the probability mass function is explicit, given in terms of the Fuss-Catalan numbers, we
can use the Zeilberger algorithm to compute recurrences for the probability of the first player winning,
for {\bf symbolic} $n$, and {\bf symbolic} $p$ (assuming that it is larger than $\frac{k}{k+1}$). Alas,
we can {\bf not} do it for symbolic $k$, since the Fuss-Catalan numbers are not bi-holonomic in both $n$ and $k$.

For the case $k=2$ we have the next proposition.

{\bf Proposition 10}: In the two player version  game, if the probability of
winning a dollar is $p$ and of losing two dollars is $1-p$ , provided $\frac{2}{3}<p<1$,
the probability of the player whose turn is to move of reaching $\geq m$ dollars first
is $(1+f(m))/2$ where $f(m)$ satisfies the sixth-order linear recurrence
$$
m \left( p-1 \right) ^{4} \left( m-5 \right) f \left( m \right) 
$$
$$
-2\, \left( p-1 \right) ^{2} \left( {m}^{2}-6\,m+6 \right) f \left( m-2 \right) 
-{p}^{2} \left( p-1 \right) ^{2} \left( 2\,{m}^{2}-13\,m+12 \right) f \left( m-3 \right) 
$$
$$
+ \left( m-3 \right)  \left( m-4 \right) f \left( m-4 \right) -{p}^{2}
 \left( 2\,{m}^{2}-15\,m+24 \right) f \left( m-5 \right) +{p}^{4} \left( m-2 \right)  \left( m-6 \right) f \left( m-6 \right) =0
\quad ,
$$
subject to the appropriate initial conditions.

For the case $k=3$  we have the next proposition.

{\bf Proposition 11}: In the two player version  game, if the probability of
winning a dollar is $p$ and of losing three dollars is $1-p$ , provided $\frac{3}{4}<p<1$,
the probability of the player whose turn is to move 
of reaching $\geq m$ dollars first is 
$(1+f(m))/2$ where $f(m)$ satisfies the eighth-order  linear recurrence
$$
m \left( p-1 \right) ^{4} \left( m-7 \right) f \left( m \right) - \left( p-1 \right) ^{2} \left( 2\,{m}^{2}-17\,m+24 \right) f \left( m-3 \right) 
$$
$$
-2\,{p}^{2} \left( p-1 \right) ^{2} \left( {m}^{2}-9\,m+12 \right) f \left( m-4 \right) + \left( m-4 \right)  \left( m-6 \right) f \left( m-6 \right) 
$$
$$
-{p}^{2} \left( 2\,{m}^{2}-21\,m+48 \right) f \left( m-7 \right) +{p}^{4} \left( m-3 \right)  \left( m-8 \right) f \left( m-8 \right) =0
\quad ,
$$
subject to the appropriate initial conditions.

For the case $k=4$  we have the next proposition.

{\bf Proposition 12}: In the two player version  game, if the probability of
winning a dollar is $p$ and of losing four dollars is $1-p$ , provided $\frac{4}{5}<p<1$,
the probability of the player whose turn is to move 
of reaching $\geq m$ dollars first is 
$(1+f(m))/2$ where $f(m)$ satisfies the tenth-order  linear recurrence
$$
m \left( p-1 \right) ^{4} \left( m-9 \right) f \left( m \right) -2\, \left( p-1 \right) ^{2} \left( {m}^{2}-11\,m+20 \right) f \left( m-4 \right) 
$$
$$
-{p}^{2} \left( p-1 \right) ^{2} \left( 2\,{m}^{2}-23\,m+40 \right) f \left( m-5 \right) + \left( m-5 \right)  \left( m-8 \right) f \left( m-8 \right) 
$$
$$
-{p}^{2} \left( 2\,{m}^{2}-27\,m+80 \right) f \left( m-9 \right) +{p}^{4} \left( m-4 \right)  \left( m-10 \right) f \left( m-10 \right) =0
\quad ,
$$
subject to the appropriate initial conditions.

For the case $k=5$ we have the next proposition.

{\bf Proposition 13}: In the two player version  game, if the probability of
winning a dollar is $p$ and of losing five dollars is $1-p$ , provided $\frac{5}{6}<p<1$,
the probability of the player whose turn is to move 
of reaching $\geq m$ dollars first is $(1+f(m))/2$ where $f(m)$ satisfies the $12^{th}$-order linear recurrence
$$
m \left( p-1 \right) ^{4} \left( m-11 \right) f \left( m \right) - \left( p-1 \right) ^{2} \left( 2\,{m}^{2}-27\,m+60 \right) f \left( m-5 \right) 
$$
$$
-2\,{p}^{2} \left( p-1 \right) ^{2} \left( {m}^{2}-14\,m+30 \right) f \left( m-6 \right) + \left( m-6 \right)  \left( m-10 \right) f \left( m-10 \right) 
$$
$$
-{p}^{2} \left( 2\,{m}^{2}-33\,m+120 \right) f \left( m-11 \right) +{p}^{4} \left( m-5 \right)  \left( m-12 \right) f \left( m-12 \right) =0
\quad ,
$$
subject to the appropriate initial conditions.

{\bf  Winning k dollars or losing one dollar}

This case is more complicated than the previous one, and we will have to treat one $k$ at a time even for
the expectation. Also, we only consider the case of reaching at least one dollar for the first time, rather
than the more general case of reaching $n$ dollars for the first time.

Now our {\bf alphabet} is $\{ \, k \, , \, -1 \, \}$ and, in terms of lattice paths, the atomic steps
are $(1,k)$ and $(1,-1)$.

Since the last step of such a path must be $(1,k)$ it can terminate at $y=k$, or $y=k-1$, \dots, $y=1$,
so we are forced to consider, in addition to $U_{0,0}$ the set of paths that start at $y=0$ and end at $y=0$ and never go above the $x$-axis, 
also $U_{0,1}$ the set of paths that start at $y=0$ and end at $y=-1$ and never go above the $x$-axis,  all the way to
$U_{0,(k-1)}$, the set of paths that start at $y=0$ and end at $y=-(k-1)$ and never go above the $x$-axis.

Such a word looks like
$$
U_{0,0} \, k \vee U_{0,1}\, k \vee \quad \dots \quad \vee U_{0,(k-1)}\, k \quad .
$$
Let $U:=U_{0,0}$. Then the weight-enumerator of $U$ is $F(z_k z_{-1}^k)$ where $F(X)$ is as above, the solution of
$$
F(X) \, = \, 1 \,+ \,X F(X)^{k+1} \quad.
$$

It can be seen that $U_{0,r}=(\1\, U_{0,0})^r$, hence its weight-enumerator is $(z_{-1} \, F(X))^r$.

Substituting for $z_{-1} \,= \, p\,t$ and $z_k \,= \, (1-p)\,t$, we get the following human-generated proposition.

{\bf Proposition 14}: 
Suppose  that at each round, you lose one dollar with probability $p$ and win $k$ dollars with probability $1-p$, and you quit 
as soon as you reach at least 1 dollar. If $0<p<\frac{k}{k+1}$ then, of course, sooner or later, you will reach your goal. 
Let $g(t)$, be the formal power series, in $t$, satisfying the algebraic equation
$$
g(t)\,- \, 1 \, - \, p^k\,(1-p)\,t^{k+1}\,g(t)^{k+1} \, = \, 0 \quad .
$$
The probability generating function, let's call it $f(t)$, for the number of rounds until having a positive capital is
$$
f(t) \, = \, (1-p)\,t\,g(t)\,\sum_{i=0}^{k-1} \,(p\,t\,g(t))^i \quad.
$$

If you will apply implicit differentiation to the defining equation of $g(t)$, and then express $f'(t)$ in terms of $g(t)$ and $g'(t)$ and
then plug-in $t=1$, you will get $0/0$. It turns out that the expressions for the expectation, variance, and higher moments
are no longer rational functions of $p$, but are roots of {\bf algebraic} equations. The reason is that when $t=1$,
$1$ is a double (or higher-order) root of the defining equation for the probability itself $f(1)=1$.

Since Maple knows how to differentiate, both explicitly and implicitly,
our beloved computer can handle it all automatically, and get explicit algebraic equation for symbolic $p$, or
specific algebraic numbers for specific $p<\frac{k}{k+1}$, alas only for {\bf one $k$ at a time}.

We have the following computer-generated proposition for the case $k=2$, i.e. for the gambling options $\{-1,2\}$, with
$Pr(-1)=p$ and $Pr(2)=1-p$.

{\bf Proposition 15}: Let $X$ be the random variable `number of rounds until you reach positive capital' if
you start at $0$, and  at each round, you lose $1$ dollar with probability $p$ and win $2$ dollars with probability $1-p$.
Assume that $p<\frac{2}{3}$.

The expectation is given by
$$
E[X] \, = \, 
{\frac {3\,p+\sqrt {\left( 3\,p+1 \right)  \left( 1-p \right) }-1}{2 \, p \left( 2 \,- \, 3\,p \right) }}
$$
For the variance, and third through the sixth moment, see

{\tt http://sites.math.rutgers.edu/\~zeilberg/tokhniot/oGenPileGames3.txt} \quad .

Note that for the most interesting case, $p=\frac{1}{2}$, the expectation is the {\bf beautiful number} $1+\sqrt{5}$
(twice the golden ratio). This is so nice that we will single it out.

{\bf Beautiful Corollary}: If a one-dimensional random walker starts at $0$ and moves {\bf one step back} with probability $\frac{1}{2}$
and {\bf two steps forward} with probability $\frac{1}{2}$ and keeps going until he is at a location $\geq 1$ for the first time,
the expected number of steps that he takes is twice the Golden Ratio, i.e. $1+\sqrt{5}$.

For $k \geq 3$ and {\it symbolic} $p$, things get too complicated to reproduce here, so  let's just mention the 
expectations for a few cases for the most interesting case, $p=\frac{1}{2}$.

$k=3$: The expected duration of a random walk with $Pr(-1)=Pr(3)=\frac{1}{2}$ until reaching a location $\geq 1$ for the first time is
the positive root of 
$$
{x}^{3}-4\,x-4=0 \quad ,
$$
that equals  $2.382975767906237494\dots$ .

$k=4$: The expected duration of a random walk with $Pr(-1)=Pr(4)=\frac{1}{2}$ until reaching a location $\geq 1$ for the first time is
the positive root of 
$$
3\,{x}^{4}+4\,{x}^{3}-8\,{x}^{2}-24\,x-16 \, = \, 0 \quad ,
$$
that equals  $2.1561901553356811691\dots$ .

$k=5$: The expected duration of a random walk with $Pr(-1)=Pr(5)=\frac{1}{2}$ until reaching a location $\geq 1$ for the first time is
the positive root of 
$$
2\,{x}^{5}+5\,{x}^{4}-20\,{x}^{2}-32\,x-16 \, = \, 0 \quad ,
$$
that equals  $2.07050432323944926\dots$ . 

$k=6$: The expected duration of a random walk with $Pr(-1)=Pr(6)=\frac{1}{2}$ until reaching a location $\geq 1$ for the first time is
the positive root of 
$$
5\,{x}^{6}+18\,{x}^{5}+20\,{x}^{4}-40\,{x}^{3}-144\,{x}^{2}-160\,x-64 \, = \, 0 \quad ,
$$
that equals  $2.0333823565252879532\dots$ . 

$k=7$: The expected duration of a random walk with $Pr(-1)=Pr(7)=\frac{1}{2}$ until reaching a location $\geq 1$ for the first time is
the positive root of 
$$
3\,{x}^{7}+14\,{x}^{6}+28\,{x}^{5}-112\,{x}^{3}-224\,{x}^{2}-192\,x-64 \, = \, 0 \quad,
$$
that equals  $2.0162018012796575781\dots$ . 

$k=8$: The expected duration of a random walk with $Pr(-1)=Pr(8)=\frac{1}{2}$ until reaching a location $\geq 1$ for the first time is
the positive root of 
$$
7\,{x}^{8}+40\,{x}^{7}+112\,{x}^{6}+112\,{x}^{5}-224\,{x}^{4}-896\,{x}^{3}-1280\,{x}^{2}-896\,x-256 \, = \, 0 \quad ,
$$
that equals  $2.00796926912597191\dots$ .

{\bf Chapter 4: Games of Pure Chance  Generated by Gambler's Ruin with Unlimited Credit: The General case}

We will now consider the general case where there is an {\bf arbitrary} set of 
non-zero integers, and an {\bf arbitrary} probability distribution on them, that we will call the {\bf die} (or {\it spinner}),
and at each round, the random walker walks (forward or backwards, as the case may be) according to the outcome of the die.
He starts at $0$, and the game ends as soon as he reaches a positive location, i.e. as soon as its location is $\geq 1$.
We will later treat the more general case where the goal is to reach a location that is $\geq m$, for any positive integer $m$.

The engine driving our algorithms is the powerful {\bf Buchberger algorithm}, that finds {\it Gr\"obner bases}, and
that is implemented in Maple and all the other major computer algebra systems.

It is convenient to separate the set of allowed steps into the set of positive steps, that we will call $U$, and
the set if negative steps, $-D$, so $D$ is a set of positive integers. For example if the set of allowed steps is
$\{-2,-5,1,3,4\}$, then $U=\{1,3,4\}$ and $D=\{2,5\}$.

Let us now state precisely the input and the desired output for our algorithms.

{\bf Main Algorithm}

{\bf Input}: 

Two sets of positive integers $D$ and $U$, corresponding to allowed steps 
$-d$ (where $d \in D$) and $u$ (where $u \in U$) in the 1D lattice, or equivalently,
$(1,-d)$ ($d \in D)$  and $(1,u)$, $u \in U$, on the two-dimensional lattice,
and an assignment of probabilities $\{p_d: d \in D\}$,
$\{p_u: u \in U\}$, such that $\sum_{ d \in D} \, p_d \, + \, \sum_{ u \in U} \, p_u \,= \, 1$ with the meaning that the
random walker walks $d$ units backward if the die landed on $d \in D$ and moves $u$ units forward if it landed
on $u \in U$. The walker starts at location $0$ and ends as soon as he reaches a strictly positive location.
In addition,  we input two {\bf symbol}s (variables), $t$ and $f$.

{\bf Output}: A polynomial $P(f,t)$ of two variables, such that
$$
P(f(t),t) \, \equiv \, 0 \quad ,
$$
holds, where $f(t)$ is the probability generating function of the random variable: `number of rounds until reaching a 
strictly positive location for the first time', obeying the above random walk.

Our algorithm guarantees that such a polynomial $P(f,t)$ {\bf always} exists.

As in Chapter 3, we will first do the corresponding {\it enumerative combinatorics} version, and later use it to our
probability purposes. We will use the powerful algorithm described in Bryan Ek's brilliant PhD thesis [Ek1], and also covered in [Ek2].

We will sometimes think of the `gambling history' listing the outcomes, getting a {\bf dynamic word} in the alphabet $U \cup (-D)$,
i.e. a 1D path, but sometimes as a static entity, its graph where $d \in D$ corresponds to the down
step $(1,-d)$ and $u \in U$  corresponds to an up-step $(1,u)$. So our problem is equivalent to counting such graphs
whose atomic steps are as above, that start at the origin, and except for the {\bf endpoint} that {\bf must} be above the $x$-axis,
is {\bf weakly below} the $x$-axis.

As before for any word $w=w_1 \dots w_n$, where $w_i$ are integers, let $Weight(w)=\prod_{i=1}^n z_{w_i}$. For example 
$Weight(1,2,-1,-3)\,=\,z_1\, z_2 \, z_{-1} \,  z_{-3}$. For any set of words $S$, its {\it weight-enumerator} is the sum of
weights of all its members. If $S$ is infinite (as is the case here) it is a {\bf formal power series} in the set of variables

$\{z_{-d} \, ; \, d \in D\} \, \cup  \{z_{u} \, ; \, u \in U\}$.

{\bf Bryan Ek's Algorithm for the Enumeration Problem}

The algorithm described in [Ek1][Ek2] does the following.

{\bf Input}: Finite sets of positive integers $D$ and $U$. This gives rise to the {\bf alphabet}
$U \cup \{-d \, : \, d \in D\}$.

{\bf Output}: A polynomial $P(f; \{z_{u}, z_{-d} \})$ of $1+|U|+|D|$, variables such that
$$
P(f (\{z_{u}, z_{-d}\} ) \,; \{ z_u, z_{-d} \} ) \, \equiv \, 0 \quad ,
$$
holds, where $f(\{z_u, z_{-d}\})$ is the weight-enumerator of all words in the alphabet $S$ whose sum is $0$ and all whose
partial sums are non-positive.

The algorithm guarantees that such a polynomial  $P(f; \{z_{u}, z_{-d} \} )$  always exists.

We use the same approach as in Chapter 3, but now we need the computer to `do the thinking', and we humans do the
`meta-thinking', teaching it how to do the `research'.

Let's abbreviate our desired weight-enumerator  $f(\{z_u, z_{-d}\})$ by $W_{0,0}$, and let
$P_{00}$ be the actual set of paths weight-enumerated by it. In other words, the set of paths starting and
ending on the $x$-axis, where each step is either $(1,u)$ ($u \in U$) or $(1,-d)$, ($d \in D$) and 
that lie {\bf weakly} below the $x$-axis.

As we will soon see, we will be forced to introduce more general quantities.
Let $W_{a,b}$ be the weight-enumerator of the set of
paths $P_{a,b}$, that start at the horizontal line $y=-a$, end at the horizontal line $y=-b$, and always stay weakly-below the
$x$-axis.

{\bf Setting up a system of Non-Linear Equations}

{\bf The case $(a,b)=(0,0)$}

Let's look at an arbitrary member, $w$, of $P_{0,0}$. It may be the {\bf empty path}, but otherwise,
let $w_1$ be the {\bf longest} prefix whose sum is $0$, then we can write
$$
w \, = \, w_1 \, w_2 \quad ,
$$
where $w_1 \in W_{0,0}$, and $w_2$ is also in $W_{0,0}$ but with the additional property that except for the endpoints,
lies {\bf strictly} below the $x$-axis.
Let's call this subset $\overline{W_{0,0}}$. Obviously, the first step of $w_2$ must be a down step, $-d$, for some $d  \in D$, and
the last step must be an up-step, $u$, for some $u \in U$. For such a path (alias word), we can write ( note that $\d=-d)$
$$
w_2 \,= \, \d \, w_3 \, u \quad,
$$
where $w_3$ is a path that starts at the horizontal line $y=-d$ and ends at the horizontal $y=-u$ but that
is {\bf strictly} below the $x-axis$. Such paths are `isomorphic' to paths that start
at $y=-(d-1)$ and end at $y=-(u-1)$ and stay weakly below the $x$-axis, in other words paths that belong to $W_{d-1,u-1}$.

So our desired quantity, $W_{0,0}$, satisfies the {\it one} non-linear equation
$$
W_{0,0} \, = \, 1 + \sum_{d \in D} \sum_{u \in U} \, z_{-d} \, W_{d-1,u-1} \, z_u \quad .
$$

Alas, now we have to handle all the `{\it uninvited guests}', the $W_{a,b}$ with $(a,b) \neq (0,0)$ that showed up.

We already handled the case $(a,b)=(0,0)$, we have to address three more cases.

{\bf The case $a>0$ and $b>0$}

If such a path, $w$, is {\bf strictly} below the $x$-axis then it is  `isomorphic' to a member of $P_{a-1 , b-1}$. 
Otherwise, sooner or later, it would meet the $x$-axis for the {\bf first time}. Let $w_1$ be the sub-path leading to that event.

We can write
$$
w \, = \, w_1 \, w_2 \quad,
$$
where $w_1$ is a path from $y=-a$ to $y=0$ that, except for the last point, lies {\bf strictly} below the $x$-axis,
let's call that set $\overline{W_{-a,0}}$.

On the other hand $w_2$ is a member of $W_{0,b}$. 
Conversely, every two such paths $w_1 \in \overline{P_{a,0}}$ and $w_2 \in P_{0,b}$,
when joined is a member of $W_{a,b}$ that touches the $x$-axis. Every path in  $\overline{W_{-a,0}}$ must obviously end with $u \in U$,
and the path obtained by removing the last step belongs to  $\overline{W_{-a,-u}}$, that is `isomorphic' to
$W_{a-1,u-1}$. We thus have the equation
$$
W_{a,b} \, = \, W_{a-1,b-1} \, + \, \left ( \sum_{u \in U} \, \, W_{a-1,u-1} \, z_u \, \right ) \, W_{0,b} \quad .
$$

{\bf The case $a>0$ and $b=0$}

The above discussion is also applicable to the case $b=0$, except that the first term on the right, $W_{a-1,b-1}$, disappears . So we have
$$
W_{a,0} \, =  \, \left ( \sum_{u \in U} \, \, W_{a-1,u-1} \, z_u \, \right ) \, W_{0,0} \quad .
$$

{\bf The case $a=0$ and $b>0$}

We have $P_{0,b}=P_{00}\, \overline{P_{0,b}}$, so similarly
$$
W_{0,b} \, =  \,  W_{0,0} \left ( \sum_{d \in D} \, z_{-d}  W_{d-1,b-1}  \right )  \quad .
$$

{\bf Symbolic Dynamical Programming}

Since our primary interest is , for now, $W_{0,0}$, the other quantities $W_{a,b}$ with $(a,b) \neq (0,0)$ are only
auxiliary unknowns, that we are not interested in for their own sake, but that would hopefully enable
us to find $W_{0,0}$.

We start out with the equation for $W_{0,0}$ that introduces $|D||U|$ new quantities, 

$$
\{ W_{d-1,u-1} \, : \, d \in D \, ,  \, u \in U \} \quad .
$$
 
For each new equation that we set-up, we may get brand new quantities, $W_{a,b}$, not yet encountered, but also some of which
that already showed up before. For each new quantity, we set up a new equation. A priori, it is conceivable that
we would have {\it infinite regress}, getting an infinite set of non-linear equations for an infinite
set of unknowns. Luckily, this does not happen! Sooner or later there are no more new `uninvited guests',
and we are left with a {\bf finite} set of non-linear (in fact quadratic) equations with the {\bf same} number of
unknowns, enabling us by {\bf elimination} (using, in our case the Buchberger algorithm in Maple)
to get {\bf one} (usually, very complicated!) equation in the one unknown, $W_{0,0}$. We can do the same thing
for any of the other $W_{a,b}$ and for that matter any linear combination of  $W_{a,b}$, calling that
linear combination $Z$, introducing one more equation and one more unknown, and {\bf eliminating} $Z$, getting
an algebraic equation satisfied by $Z$.

{\bf Straight Enumeration}

Suppose that you are interested in the actual {\it enumerating sequence}, i.e. given a set of non-negative integers, $S$,
you want to have an ``explicit" expression for $a(n)$ the number of 1D walks, starting at $0$ using the steps of $S$,
ending at $0$, and always staying weakly to the left of the origin. Equivalently, given a set of integers $S$,
our `alphabet', $a(n)$ is the number of words of length $n$ whose sum is $0$ and all whose partial sums are non-positive.

By {\bf specializing} $z_u=t$, ($u \in U)$ ;  $z_d=t$, ($d \in D)$ we have our next conceptual theorem.

{\bf Theorem 5}: For an {\bf arbitrary} set of integers $S$, let $a_S(n)$ be the number of sequences of length $n$, whose
entries are drawn from $S$ with the property that the sum is $0$ and all its partial sums are non-positive. Then the ordinary generating function,
in the variable $t$,
$$
f(t) \, := \, \sum_{n=0}^{\infty} a_S(n) t^n \quad,
$$
is an algebraic formal power series, i.e. there exists a two-variable polynomial $P$, with {\bf integer coefficients}, such that
$$
P(\, f(t) \, , \, t \, ) \, \equiv \, 0 \quad .
$$
Furthermore, there exists an algorithm for finding this polynomial $P(f,t)$.

As a corollary, we know that $f(t)$ is $D-finite$, and hence $a(n)$ is $P$-recursive, and we have

{\bf Theorem 5'}: For an {\bf arbitrary} set of integers $S$, let $a_S(n)$ be the number of sequences of length $n$, whose
entries are drawn from $S$ with the property the sum is $0$, and that all the partial sums are non-positive.  The sequence $a_S(n)$ is
{\bf P-recursive}, i.e. there exists a positive integer $L$ and polynomials $p_i(n)$ in $n$, such that
$$
\sum_{i=0}^{L} \, p_i(n) a_S(n-i) \, \equiv \, 0 \quad .
$$
Furthermore, there exists an algorithm for finding this linear recurrence.

{\bf Rigorous Experimental Mathematics}

Now that we have the {\bf theoretical guarantee}  that the polynomial $P(f(t),t)$ and the recurrence {\bf exists},
it may be more efficient, to crank out, using (the usual, not symbolic) {\bf dynamical programming}, sufficiently
many terms and then fit them with a recurrence.  
The Maple commands  {\tt listtoalgeq(l,y(x))} and {\tt listtorec(l,a(n))}  do just that.
These two useful commands are part of the versatile Maple package {\tt gfun} written
by Bruno Salvy and Paul Zimmerman [SZ].

{\bf The Weight-Enumerator of 1D walks that start at $0$ end at strictly positive location but
otherwise stay weakly to the left of $0$}

Equivalently, in the two-dimensional version, our walks start at the origin, end above the $x$-axis,
and except for the end-point are weakly below the $x$-axis.

Since every such walk must obviously end with an up-step, $u \in U$, and the endpoint could be either at
$y=1 , y=2, \dots, y=u$, 
the desired weight-enumerator, let's call it $Z$, using the $W_{a,b}$ above is
$$
Z \, = \, \sum_{u \in U} \, \left ( \sum_{u'=0}^{u-1} \, W_{0,u'}  \right )  \, z_u\quad.
$$
This is our quantity $Z$ mentioned above, and thanks for the Buchberger algorithm, we can eliminate everything except $Z$,
and get a {\bf single} polynomial equation in $Z$.

{\bf From Enumeration to Probability}

Having gotten a polynomial equation satisfied by the formal power series in the set of $|D|+|U|$ variables
$\{ z_u \, : \, u \in U\} \, \cup \,  \{ z_{-d} : \, d \in D\}$ that enumerates the above set of words,
we {\bf plug-in} $z_u\,=\,p_u \,t$, $z_{-d}\,=\,p_d \,t$, getting an equation of the form
$$
P(f(t),t) \, \equiv \, 0 \quad ,
$$
satisfied by the probability generating function, $f(t)$, for the random variable `duration until reaching a positive amount for
the first time' if you start at $0$ dollars, at each time step (round) you win $u$ dollars with probability $p_u$ , if $u \in U$,
and lose $d$ dollars with probability $p_d$, if $d \in D$.

If the expected gain of a single step
$$
\sum_{u \in U} \, p_u \,u \, - \,\sum_{d \in D} \, p_d \,d \ \quad,
$$
is positive, then sooner or latter the game ends.

So far $f(t)$ was a {\bf formal power series}, but when you plug-in $t=1$, you get a (numerical) {\bf convergent} power series that
must sum to $1$, so $f(1)=1$. It follows, that $1$ is one of the roots of
the {\bf numerical} equation , in $f(1)$.
$$
P(f(1),1) \, = \, 0 \quad ,
$$

This implies our next {\bf conceptual} result, Theorem 6.

{\bf Theorem 6}: Let $\P$ be any `die' (with {\bf any}  number of faces, {\bf any} loading, and {\bf any} assignments of non-zero integers to its faces),
where, at each step, you win or lose according to the outcome.
Assume that the expected gain of a single round is positive.
Let $X$ be the random variable `number of rounds until reaching a positive amount for the first time'. 

The probability generating function of $X$,
$$
f(t) \, = \, \sum_{k=0}^{\infty} Prob(X=k) \, t^k \quad,
$$
is an {\bf algebraic} formal power series, in other words, there exists a two-variable polynomial, $P(f,t)$, such that
$$
P(f(t), t) \equiv 0 \quad .
$$
Furthermore, there is an algorithm for finding  the two-variable polynomial $P(f,t)$.

As a corollary, we know that $f(t)$ is $D-finite$, and hence $Prob(X=k)$ is $P$-recursive, in $k$.

{\bf Theorem 6'}: For an {\bf arbitrary} `die' as in Theorem 6, and $X$ defined there, the sequence $a(k)=Prob(X=k)$
is {\bf P-recursive}, i.e. there exists a positive integer $L$ and polynomials $p_i(k)$ in $k$, such that
$$
\sum_{i=0}^{L} \, p_i(k) \, a(k-i) \, \equiv \, 0 \quad .
$$
Furthermore, there exists an algorithm for finding this linear recurrence.

What about expectation?, using {\it implicit differentiation}, we get
$$
P_f(f,t) \cdot f'(t) + P_t(f,t) \equiv 0 \quad .
$$
Eliminating $f=f(t)$, from the two equations with three variables $\{f,f',t\}$ ($f'$ is short for $f'(t)$)
$$
P_f(f,t) \cdot f' + P_t(f,t) \,= \,0 \quad, \quad P(f,t) \, = \, 0 \quad .
$$
we (or rather our computers) get a polynomial equation of the form $Q(f'(t),t) \, = \, 0$,
and plugging-in $t=1$, ($f'(1)$ is a numerical convergent series),
we get that the expected duration is one of the roots of the {\bf numerical} equation
$$
Q(f'(1),1) \, = \, 0  \quad .
$$ 
By repeated implicit differentiation, and elimination of $(t\frac{d}{dt})^2 f(t)$, 
and then $(t\frac{d}{dt})^3 f(t)$, etc., we get algebraic equations for as many moments as desired, and hence for the variance,
and higher moments about the mean. 
All this is implemented in procedure {\tt Momk(N,P,fk,k)} in the Maple package {\tt VGPileGames.txt}.
Readers who wish to see more details are more than welcome to examine the Maple source code.

This brings us to the next `conceptual' result.

{\bf Theorem 7}: Consider {\bf any} finite set of non-zero integers, and any probability distribution on them with
positive expectation, where, at each step, you win or lose according to the outcome.
Assume that the expected gain of a single round is positive.
Let $X$ be the random variable `number of rounds until reaching a positive amount for the first time.'
Then the expectation, variance, and any higher moments of $X$ are {\bf algebraic numbers}, whose
minimal polynomial can be explicitly computed.

If the expected gain of a single round, $\sum_{u \in U} \, p_u \,u \, - \,\sum_{d \in D} \, p_d \,d$, is $0$, then
the expectation and higher moments are infinite. If it is negative, then the probability of exiting
with a  positive amount is less than $1$, and $f(1)$ is one of the roots of the numerical equation $P(f(1),1)=0$, that
is an explicit algebraic number. Then one has to talk about the `conditional duration', and replace
$f(t)$ by $f(t)/f(1)$, and Theorems 6 and 7 are still valid.

We will only present here one case. Quite a few similar propositions can be found in the web-page of this article, and readers
are welcome to  generate many more using the command 

{\tt Paper(N,P,k,K1,K2,f,eps);}

in the Maple package {\tt VGPileGames.txt} also available from there.

{\bf Proposition 16}: Consider a 1D random walk with a set of steps $\{-1,-2,1,2\}$ where 
$Pr(-1)=\frac{1}{4}$, $Pr(-2)=\frac{1}{8}$, $Pr(1)=\frac{1}{4}$, $Pr(2)=\frac{3}{8}$,
that starts at $0$. Let $X$ be the random variable:

`number of steps until reaching a strictly positive location for the first time'.

The probability generating function of $X$
$$
f(t) \, = \, \sum_{k=0}^{\infty} Prob(X=k) \, t^k \quad,
$$
is a formal power series that satisfies the algebraic equation
$$
{t}^{3}{f}^{6}+ \left( 6\,{t}^{3}-8\,{t}^{2} \right) {f}^{5}+ \left( 
19\,{t}^{3}-48\,{t}^{2} \right) {f}^{4}+ \left( 84\,{t}^{3}-80\,{t}^{2
}+128\,t \right) {f}^{3}
$$
$$
+ \left( 71\,{t}^{3}-608\,{t}^{2}+320\,t
 \right) {f}^{2}+ \left( 262\,{t}^{3}-360\,{t}^{2}+768\,t-512 \right) 
f+69\,{t}^{3}-432\,{t}^{2}+320\,t
\, = \, 0 \quad.
$$

The expectation, $f_1$, is one of the roots of the cubic equation
$$
{f_{{1}}}^{3}-12\,{f_{{1}}}^{2}+16\,f_{{1}}+32 \, = \, 0 \quad,
$$
whose floating-point rendition is  $2.9653919099833889\dots$.

The second moment, $f_2$, is one of the roots of the cubic equation
$$
101\,{f_{{2}}}^{3}-14140\,{f_{{2}}}^{2}+367216\,f_{{2}}-273824 \, = \, 0 \quad .
$$
whose floating-point rendition is  $33.31799734943726426\dots$. It follows that the variance is
$24.5244481696423327\dots$, and hence the standard-deviation is
$4.9522164905870515\dots.$

{\bf The two-player version}

Suppose two players take turns rolling the $\P$ die, and the one who is the first to reach a positive
amount is declared the winner.
Recall that the probability of the first player winning the game is
$(1+s)/2$, where
$$
s= \sum_{n=1}^{\infty} Prob(X=k)^2 \quad ,
$$
it follows from Theorem 6' that $s$ is a {\bf holonomic constant}.

{\bf Playing until reaching at least $m$ dollars for the first time}

If the set of up-steps, $U$, consists of only one dollar, i.e. if $U=\{1\}$, then the probability generating function
for the random variable ``number of rounds it takes until reaching at least $m$ dollars for the first time''
is simply $f(t)^m$, and everything goes through, and the probability of the first player winning is
holonomic in $m$. In the more general case, it is also true, but a bit more complicated, and we did
not implemented it yet. The probability generating function for the  `first time of reaching $\geq m$' case can be
shown to satisfy a linear recurrence equation  in $m$ whose coefficients are  what we called
$\{W_{a,b}\}$ above. After differentiating with respect to $t$,
we can get recurrences for the expectation, and higher moments.

This brings us to the next theorem that we state without proof, and is not yet implemented in general.

{\bf Theorem 8}: Consider {\bf any} finite set of non-zero integers, and any probability distribution on them with
positive expectation, where, at each step, you win or lose according to the outcome.
Assume that the expected gain of a single round is positive. For any positive integer $m$,
let $X_m$ be the random variable `number of rounds until reaching an amount that is $\geq m$ for the first time.'
Then the probability generating function of $X_m$, let's call it $f_m(t)$ is an (constant) algebraic formal power series,
and $f_m(t)$ satisfies a linear recurrence in $m$ with coefficients that are algebraic formal power series.

Furthermore, $E[X_m]$ are {\bf algebraic numbers}, that satisfy a linear recurrence equation
with {\bf constant} coefficients (but the constants featuring in the linear recurrence are, in general,
{\it algebraic} numbers).

{\bf Keeping it Simple: Numerics Driven by Symbolics}

For quite a few `dice', our computers were able to find the {\bf exact answer} for the question 

{\it
What is the probability generating function for the random variable 
`Number of rounds until reaching a positive amount for the first time \quad .
}

Of course, except for the Catalan case, it is not {\bf fully} explicit, but it is {\bf as explicit as it gets}, 
the {\bf exact} polynomial equation
$$
P(f(t),t) \, = \, 0 \quad,
$$
satisfied by it, and that, in turn, enables us to find the {\bf exact} values of the expectation, variance and higher moments,
in terms of {\bf explicit} algebraic numbers, i.e. numbers given by their minimal equation with integer coefficients.

But if the `die', $\P$, gets larger, these algorithms are mainly of {\it theoretical} interest, i.e. for
{\bf conceptual computation}. To actually get {\bf answers}, {\bf very fast}, we recommend using
the following simple-minded {\it symbolic-numeric} algorithm.

We can also do {\bf simulation}, but these are very inexact. They are only useful (for our current project)
as {\bf sanity checks}, to make sure that we did not mess up.

Let $h(x)$ be the {\bf probability generating function} of our die
$$
h(x)\,=\, \sum_{u \in U} p_u \, x^u \, + \,  \sum_{d \in D} p_d \, x^{-d} \quad .
$$
For example, for the fair Catalan case $h(x)=\frac{1}{2}(x+x^{-1})$, for the Fuss-Catalan case considered in Chapter 3,
with `one step forward, $k$ steps backwards', $h(x)=px+(1-p)x^{-k}$, and for the more difficult case
`one step backwards, $k$ steps forward' case, we have  $h(x)=px^{-1}+(1-p)x^k$.

For any Laurent polynomial, define the operator: `the positive part' as follows:
$$
G[\sum_{i=c}^{d} a_i x^i) ] \, = \,\sum_{i=1}^{d} a_i x^i \quad .
$$
For example, 
$$
G[\frac{1}{10}x^{-3} + \frac{1}{20}x^{-2}+ \frac{7}{20}x^{-}+ \frac{1}{20}+\frac{1}{4}x+ \frac{1}{4}x^2]
=\frac{1}{4}x+\frac{1}{4}x^2 \quad .
$$

{\bf The Symbolic-Numeric Algorithm to compute the first $K$ terms of the probability generating function of our `duration of the game' random variable }

{\bf Input}

$\bullet$  A die, $\P$, whose probability generating function is the Laurent polynomial $h(x)$.

$\bullet$ A positive integer $K$.

{\bf Output}

The first $K$ terms in the Maclaurin expansion of the probability generating function of the random variable:
`{\it number of rounds until the player reaches a strictly positive amount for the first time}',
let's call it $f_K(t)$.

{\it Initialize}: $F_0(x):=1$, $f_0(t)=0$.

For $i$ from 1 to K do

$$
A(x):= \, F_{i-1}(x)\, h(x) \quad ,
$$
$$
F_{i}(x)\, = \, A(x) \, - \, G[A(x)] \quad ,
$$
$$
f_i(t)=f_{i-1}(t) + G[A(x)]|_{x=1} \, t^i \quad .
$$

Intuitively, $F_i(x)$ describes the scenarios that still did not make it to positivity by the $i$-th step. Multiplying
by $h(x)$ is the `roll of the die', $G[A(x)]$ describes the lucky scenarios that made it by the $i$-th round, and plugging in
$x=1$, gives the probability due to all the scenarios that made it exactly at the $i$-round for the first time.

If $h'(1)>0$, then $f(1)=1$, and to see how good $f_K(t)$ approximates $f(t)$, plug-in $t=1$.
If this is very close to $1$ (usually it is!), then you are safe.

Also, {\it frankly}, you are {\bf not} immortal, and even if you are, it is good to {\bf a priori} set a
limit to the number of allowed rounds, and compute everything conditioned on finishing in $\leq K$ rounds.

The conditional expectation on finishing in $\leq K$ rounds is $f_K'(1)/f_K(1)$, 
the second moment is $(t\frac{d}{dt})^2f_K(t)|_{t=1}/f_K(1)$
and the $k$-th moment is $(t\frac{d}{dt})^kf_K(t)|_{t=1}/f_K(1)$.

These give much faster, very accurate, approximations to the desired statistical quantities of our random variable.

Similarly, we can compute, very fast, the truncated Taylor series of the random variable
`first time of having an amount $\leq m$', for any desired $m$.

This is accomplished by procedure {\tt Ngf(N,P,t,K)} in the Maple package {\tt VGPileGames.txt}.

If you want to see many examples, look at the file

{\tt http://sites.math.rutgers.edu/\~{}zeilberg/tokhniot/oVGPileGames2.txt} \quad .

Another route to get the {\bf exact} value of the expectation, variance, etc., is to
derive numerical approximations like we did, and use Maple's command {\tt identify}
or use the {\bf Inverse  Symbolic Calculator},

{\tt https://isc.carma.newcastle.edu.au/} \quad .

If  you get an algebraic number, then it is most likely the right one, since we know,
from the `conceptual part', that it {\bf is} an algebraic number.

Using this `experimental way', we discovered the following {\bf lovely} proposition.

{\bf Proposition 17}: Consider a `two steps forward one step backwards random walk' starting at $0$,
with $Prob(-1)=Prob(2)=\frac{1}{2}$. The expected number of rounds until reaching a location $\geq m$
for the {\bf first} time, equals, {\bf exactly}
$$
2\,m \, + \, (4\,-\, 2 \, \phi) \, + \,  2\,(F_{m+2} \, \phi -F_{m+3}) \quad ,
$$
where $F_m$ are the Fibonacci numbers and $\phi=\frac{1+\sqrt{5}}{2}$ is the Golden Ratio.

Note that this makes sense, since the last term $F_{m+2} \, \phi -F_{m+3}$ is exponentially small in $m$, so
this is very close to $2\, m \,+ \, 4 \,- \,2 \phi$, and since the expected gain of one round is $\frac{1}{2}$ a crude approximation
to the expected duration until owning $\geq m$ dollars should be roughly $m/\frac{1}{2}\,=\,2m$.
Also note that when $m=1$ we get $2+4-2\phi+ 2(F_3\,\phi\, - \, F_4)=6-2\phi+ 2(2\,\phi\, - \,3)\,=\, 2\phi$, in agreement
with the result established in Chapter 3.

{\bf What about a proof of Proposition 17?}

Proposition 17 was discovered {\it experimentally}, but we {\bf do} know how to prove it. Writing it up, though, will take
time and effort that we are unwilling to spend. We will be glad to furnish a proof in return to
a $\$2000$ donation to the {\it On-Line-Encyclopedia of Integer Sequences}. 

{\bf Conclusion} 

In this article, using {\it Games of pure chance} as a {\bf case study}, we preached the
value of {\bf computational diversity}. Purely numeric, numeric-symbolic, purely symbolic, and `conceptual',
as well as the {\bf simulation}, that in 
our case plays a secondary role, as a {\bf checker}. It is so easy to have bugs in your programs, or
gaps in your reasoning, so it is still reassuring that you can confirm the numbers that you got
are in the right ball-park.
We also demonstrated a novel application of the {\bf Buchberger algorithm}.

{\bf References}

[AZ] Gert Almkvist and Doron Zeilberger, {\it The Method of differentiating Under The integral sign},
J. Symbolic Computation {\bf 10} (1990), 571-591. Available from \hfill\break
{\tt http://sites.math.rutgers.edu/\~{}zeilberg/mamarim/mamarimhtml/duis.html } \quad .

[DM]  A. Dvoretzky and Th. Motzkin, {\it A problem of arrangements}, Duke Math. J. {\bf 14} (1947), 305-313.

[Ek1] Bryan Ek, ``Unimodal Polynomials and Lattice Walk Enumeration with Experimental Mathematics'', 
PhD thesis, Rutgers University, May 2018. Available from \hfill\break
{\tt http://sites.math.rutgers.edu/\~{}zeilberg/Theses/BryanEkThesis.pdf} \quad .

[Ek2]  Bryan Ek, {\tt Lattice Walk Enumeration}, 29 March, 2018, {\tt https://arxiv.org/abs/1803.10920}.

[F] William Feller, ``{\it An Introduction to Probability
Theory and Its Application}'', volume 1, three editions.
John Wiley and sons. First edition: 1950. Second edition:
1957. Third edition: 1968.

[KP] Manuel Kauers  and Peter Paule, {\it ``The Concrete Tetrahedron''}, Springer, 2011.

[LT]  Ho-Hon Leung and Thotsaporn ``Aek'' Thanatipanonda, {\it A Probabilistic Two-Pile Game},
Journal of Integer Sequences, {\bf 22\#4} (2019). \hfill\break
Also available from {\tt https://arxiv.org/abs/1903.03274} $\,\,$.

[MM] David Mccune and Lori Mccune, {\it Counting your chickens with Markov chains}, Mathematics Magazine {\bf 92} (2019), 162-172.

[PWZ] Marko Petkovsek, Herbert S. Wilf, and Doron Zeilberger, {\it ``A=B''}, A.K. Peters, 1996. Freely available from 
{\tt https://www.math.upenn.edu/\~{}wilf/AeqB.html} \quad .

[SZ] Bruno Salvy and Paul Zimmerman, {\it GFUN: a Maple package for the manipulation of generating and holonomic functions in one variable},
ACM Transactions on Mathematical Software {\bf 20}(1994), 163-177 .

[Sl] Neil J. A. Sloane, {\it The On-Line Encyclopedia of Integer Sequences}, {\tt http://www.oeis.org} .

[St] Richard Stanley, ``{\it Catalan Numbers}'', Cambridge University Press, 2015.

[T] Thotsaporn ``Aek'' Thanatipanonda, {\it A Quantitative Study on Average Number of Spins of Two-Player Dreidel},
{\tt https://arxiv.org/abs/1907.11851 } \quad .

[Z1] Doron Zeilberger, {\it  The C-finite Ansatz}, Ramanujan J. {\bf 31}(2013), 23-32.
Available on-line: \hfill\break
{\tt http://www.math.rutgers.edu/\~{}zeilberg/mamarim/mamarimhtml/cfinite.html}

[Z2] Doron Zeilberger, 
{\it  An Enquiry Concerning Human (and Computer!) [Mathematical] Understanding},
Appeared in: C.S. Calude, ed., ``Randomness \& Complexity, from Leibniz to Chaitin'' World Scientific, Singapore, 2007.
Available on-line: \hfill\break
{\tt http://www.math.rutgers.edu/\~{}zeilberg/mamarim/mamarimhtml/enquiry.html}

[Z3] Doron Zeilberger, {\it A holonomic systems approach to special function identities}, J. Computational
and Applied Mathematics {\bf 32} (1990), 321-368. Available on-line: \hfill\break
{\tt http://sites.math.rutgers.edu/\~{}zeilberg/mamarim/mamarimhtml/holonomic.html} \quad .

[Z4] Doron Zeilberger, {\it  Lagrange Inversion Without Tears (Analysis) (based on Henrici)},
The Personal Journal of Shalosh B. Ekhad and Doron Zeilberger \hfill\break
{\tt http://sites.math.rutgers.edu/\~{}zeilberg/mamarim/mamarimhtml/lag.html} \quad .

\vfill\eject

\bigskip
\hrule
\bigskip
Thotsaporn ``Aek'' Thanatipanonda, Mahidol University International College, Nakornpathom, Thailand \hfill\break
Email: {\tt thotsaporn at gmail dot com} \quad .
\bigskip
Doron Zeilberger, Department of Mathematics, Rutgers University (New Brunswick), Hill Center-Busch Campus, 110 Frelinghuysen
Rd., Piscataway, NJ 08854-8019, USA. \hfill\break
Email: {\tt DoronZeil at gmail  dot com}   \quad .
\bigskip
\hrule
\bigskip
Written: Sept. 24, 2019.
\end